\input amstex
\documentstyle{amsppt}
\magnification=\magstep1
\vcorrection{-0.8cm}
\TagsOnRight

\nologo
\def\nd{\noindent}
\def\ni{\noindent}
\def\parni{\par\ni}

\def\endemo{\qed\enddemo}

\def\sbs{\subset}
\def\sps{\supset}
\def\sbsq{\subseteq}
\def\spsq{\supseteq}

\def\a{\alpha} 
\def\b{\beta} 
\def\g{\gamma}
\def\d{\delta} 

\def\th{\theta} 
\def\l{\lambda}
\def\s{\sigma}
\def\t{\tau}
\def\u{\upsilon}

\def\thvee{\theta^\vee}
\def\o{\omega^\vee}

\def\Q{ Q^\vee} 

\def\p{\Phi}
\def\D{\Delta} 
\def\Dp{\Delta^+} 
\def\Dm{\Delta^-} 
\def\Da{\widehat\Delta}
\def\Dap{\widehat\Delta^+} 
\def\Pia{\widehat \Pi}

\def\gog{\frak g} 
\def\goh{\frak h}
\def\gok{\frak k}
 
\def\gob{\frak b} 
\def\goi{\frak i}
\def\goj{\frak j}

\def\gom{\frak m}

\def\wp{W_\a}
\def\wap{\widehat W_\a}
\def\dperp{\D_\a}
\def\dap{\Da_\a}

\def\i{\frak i}

\def\wi{w_\goi}

\def\nat{\Bbb N} 
\def\ganz{\Bbb Z} 
 
\def\real{\Bbb R}

\def\Npm{\overline N}

\def\weak{\leq}

\def\supp{\text{\it Supp}}
\def\stab{\text{\it Stab}}

\def\max{\text{\it max }}
\def\min{\text{\it min }}

\def\dim{\text{\it dim }}

\def\Waf{W_{\hskip-2pt\text{\it af}}}
\def\wab{{\Cal W}_{\hskip-2pt\text{\it ab}}}
\def\hatw{\widehat W}

\def\I{\Cal I}
\def\II{\I_{\hskip-1pt\text{\it ab}}}

\def\dab{D_{\hskip-2pt\text{\it ab}}}

\def\({\left(} 
\def\){\right)} 

\def\la{\langle}
\def\ra{\rangle}

\def\ov{\overline}

\def\y{y_{\t,\a}}

\def\wa{\check w_\a}
\def\wah{\check w_{h,\a}}

\def\adj{\text{\it Adj}}

\def\meno{\hskip.5pt\raise 2.5pt\hbox{{$\scriptscriptstyle\setminus$}}}

\topmatter 

\author  
Paola Cellini\\ Paolo Papi
\endauthor 
\abstract
This paper is devoted to a detailed study of certain remarkable posets
which form a natural partition  of all abelian ideals of a Borel
subalgebra. Our main result is a  nice uniform formula for the
dimension of maximal  ideals in these posets. We also obtain results
on $ad$-nilpotent ideals which complete the analysis started in
\cite{CP2}, \cite{CP3}.
\endabstract
\address 
Paola Cellini 
\vskip 0.pt Dipartimento di Scienze
\vskip 0.pt 
Universit\`a ``G. d'Annunzio" 
\vskip 0.pt 
Viale Pindaro 87
\vskip 0.pt 
65127 Pescara --- ITALY 
\vskip 0.pt 
e-mail:{\rm \ cellini\@sci.unich.it } 
\endaddress 

\bigskip 
\bigskip 
\address 
Paolo Papi
\vskip 0.pt 
Dipartimento di Matematica Istituto G. Castelnuovo 
\vskip 0.pt
Universit\`a di Roma ``La Sapienza" 
\vskip 0.pt 
Piazzale Aldo Moro 5 
\vskip 0.pt
00185 Roma --- ITALY 
\vskip 
0.pt e-mail:{\rm \ papi\@mat.uniroma1.it}
\endaddress

\leftheadtext 
{Paola Cellini Paolo Papi} 
\rightheadtext 
{$ad$-nilpotent ideals of  a Borel subalgebra 
III} 

\title $ad$-nilpotent ideals of  a Borel subalgebra 
III\endtitle



\endtopmatter 
\document

\heading 
\S1 Introduction
\endheading
\medskip

In this paper  we  investigate in detail certain remarkable posets of
abelian ideals $\II(\t,\a)$ of a Borel subalgebra of a simple Lie
algebra.  These posets depend upon a long root $\a$ and an element
$\t$ of the coroot lattice (subject to certain natural restrictions),
and  form a partition  of all  abelian ideals of a Borel subalgebra.
As a consequence of our analysis we may deduce  uniform formulas for
the dimensions of maximal and minimal ideals in $\II(\t,\a)$. Our work
generalizes and refines the work of Panyushev and Suter on similar
topics, using as a key tool the combinatorics of affine Weyl groups
and root systems. Our formulas also provide a uniform explanation of
classical Malcev's results on the maximal dimension of an abelian
subalgebra of a simple Lie algebra. In order to explain more precisely
our results we need to set the basic framework.
\par

Let $\gog$ be a complex simple Lie algebra of rank $n$. Let
$\gob\subset\gog$ be a fixed Borel subalgebra, with Cartan component
$\goh$; let $\Dp$ be the positive system of the root system $\D$ of
$\gog$ corresponding to the previous choice, and
$\Pi=\{\a_1,\ldots,\a_n\}$ be the associated set of simple roots. Let
$Q,\,Q^\vee$ be the root and coroot lattices, and $W$ be the Weyl
group of $\gog$. For each $\a\in \Dp$ let ${\frak g}_\a$ denote the
root space of $\gog$ relative to $\a$. Denote by $\II$ the set of
abelian ideals  of $\gob$.
\smallskip 

The class $\II$  has been introduced in \cite{K1}, but the paper that
started the recent research activity on the subject is the long
research announcement \cite{K2}. In this paper Kostant established 
deep connections between the combinatorics of abelian ideals, the
machinery of Cartan decompositions and the theory of discrete series
(for developments in this direction see also \cite{CMP}) . His work is
partly  based on  a simple but surprising result due to D. Peterson
(the so-called Peterson's $2^{rk(\gog)}$ abelian ideals theorem): the
abelian ideals may be encoded and in turn enumerated through elements
in a certain subset $\Cal W_{ab}$ of the affine Weyl group
$\widehat{W}$  of $\gog$. In \cite{CP2} we extended this encoding, say
$\goi\mapsto w_\goi$,  to the set $\I$ of ideals of $\gob$ contained
in $[\gob,\gob]$. These ideals are characterized, among the ideals of
$\gob$, by the following property: for any $x\in\goi$,
$ad(x):\gog\to\gog$ is a nilpotent endomorphism,  hence we call them 
{\it $ad$-nilpotent ideals}. Moreover we gave a characterization of
the image of the map $\goi\mapsto w_\goi$ inside $\widehat W$.
Reformulating  this characterization in terms of the affine
representation of $\widehat{W}$ allowed us to describe  a suitable
simplex $D$ in $\frak h_{\Bbb R}^*$ whose $Q^\vee$-points parametrize
$\I$ (see \cite{CP3}). In such a way we obtained a bijection between
$\I$ and certain $W$-orbits in the coroot lattice  of $\gog$, which
affords a uniform enumeration of $\I$. A short part of section  3 of 
this paper is devoted to give results on $\I$  which complete the
analysis started in our previous papers and which have been announced
in \cite{CP3}: we characterize the ideals in $\I$ having a fixed class
of nilpotence in terms of elements of $D\cap Q^\vee$ and the
generators of $\goi\in\I$ as a $\gob$-module in terms of $w_\i$.
\smallskip 

The remaining and  more extended part of the paper deals with abelian
ideals. Recall from \cite{CP2} that if $\goi\in\II$ and $\wi=t_\t v,$
where $t_\t$ is the translation by $\t\in\Q,$ and $v\in W$, then $\t$
belongs to set $X$ of $\Q$-points lying in twice the closure of the
fundamental alcove of $\widehat W$. We provide a general description
of the elements of $X$  as orthogonal sums of highest roots of nested
standard irreducible parabolic subsystems of $\D$. This result  turns
out to  be a very useful technical device  for our goals.\par We begin
our analysis decomposing $\II$ w.r.t. $X$: the map $\Cal W_{ab}\to
X,\, t_\t v\mapsto \t$ is surjective, thus we have a partition into
non-empty sets
$$\II=\bigcup_{\t\in X} \II(\t),$$ 
$\II(\t)$ being the set of ideals $\goi$ whose $\wi$ has $t_\t$ as
translational part. In \cite{K2} and \cite{CMP} it is shown that
$\II(\t)$ is in canonical bijection with $W_\tau\backslash W$, where,
if $\frak k_\tau$ is the equal rank symmetric subalgebra of $\gog$
formed by the fixed points of $Ad(\exp(\sqrt{-1}\pi\tau)$, $W_\tau$
denotes the Weyl group of $(\gok_\tau,\goh)$. In the following
discussion we list our main results.
\smallskip

We prove that, for all $\t\in X,$ $\II(\t)$ has a unique minimal
element; moreover we describe explicitly the maximal elements  of
$\II(\t)$. If $\t=0$, then $\II(\t)$ consists just of the zero ideal, 
so assume $\t\ne 0$. If $\goi$ is maximal in $\II(\t)$, then
$w_\i^{-1}(-\th+2\d)$ is a long simple root (here $\th$ is the highest
root of $\D$ and $\d$ is the fundamental imaginary root of the affine
system $\Da$ attached to $\D$). So we are led to consider, for any
$\a$ in the set $\Pi_\ell$ of long simple roots, the subposets
$$\II(\t,\a)=\{ \goi\in \II\mid w_\i=t_\t v,\,
w_\i^{-1}(-\th+2\d)=\a\}.$$ 
For any $\goi\in \II(\t)$, $w_\i^{-1}(-\th+2\d)$ belongs to $\D$ and
is a long positive root. We denote by $\Dp_\ell$ the set of long
positive roots in $\D$ and extend the definition of $\II(\t, \a)$ to
all $\a$ in $\Dp_\ell$. Thus we have
$$\II(\t)=\bigcup\limits_{\a\in\Dp_\ell} \II(\t,\a).$$
$\II(\t,\a)$ is possibly empty; we give an explicit criterion for
$\II(\t,\a)$ to be non-empty. We prove that if
$\II(\t,\a)\ne\emptyset$, then it has unique minimal and maximal
elements. More precisely, we prove that  the poset $\II(\t,\a)$ is
order isomorphic to the coset space $W((\Pi\cap\a^\perp)\setminus
P_{\t,\a})\backslash W(\Pi\cap\a^\perp)$ (see 7.6 for the precise
definition of $P_{\t,\a}$ and 2.8, 2.9 for the definition of the poset
of cosets.) This result is a generalization of \cite{S, Theorem 11}.
We explicitly describe elements $w_\t,\,y_{\t,\a},\,\check
y_{\t,\a}\in \widehat W$ such that
$$\i(w_\t)=\min\II(\t),\quad \i(w_\t y_{\t,\a})=\min\II(\t,\a),\quad
\i(w_\t \check y_{\t,\a})=\max\II(\t,\a).$$
The ideals $\max\II(\t,\a)$ for all simple roots $\a$ such that
$\II(\t,\a)\ne\emptyset$ are exactly the maximal elements of
$\II(\t)$.
\smallskip

For $\a\in\Dp_\ell$ set $X_\a=\{\t\in X\mid \II(\t,\a)\ne
\emptyset\}$.  Then $X_\a$ is  chain in $\Q$, in the standard partial
order, and its minimal element is $\thvee$. Assume
$X_\a=\{\t_1,\ldots,\t_{k(\a)}\},\, \thvee=\t_1<\ldots<\t_{k(\a)}$,
and denote by $\gom_h(\a)$ the maximal ideal in $\II(\t_h,\a)$. Then a
careful analysis of the elements of $\widehat W$ corresponding to the
ideals $\gom_h(\a)$ allows us to prove that
$$\gom_1(\a)\sbs\ldots\ldots\sbs \gom_{k(\a)}(\a).$$ 
The proof of these inclusion relations gives  a clear explanation of
the ``shell-like" structure of the elements of $\wab$ noted in
\cite{S, \S5}. Moreover it gives as a by-product a  case-free proof of
the following nice result, due to Suter \cite{S, Corollary 13}, which
holds for the maximal abelian ideals
$\gom(\a)=\gom_{k(\a)}(\a),\,\a\in\Pi_\ell$ of $\gob$:
$$\dim(\gom(\a))=g-1+\frac{1}{2}\left(|\Da_{\a}|-|\D_{\a}|\right)$$
where $\Da_{\a},\,\D_{\a}$ are the (finite) subsystems generated by
the simple roots orthogonal to $\a$ in $\Da,\,\D$ respectively, and
$g$ is the dual Coxeter number of $\D$.
\par

In particular we obtain a uniform version of  Malcev's formulas for
the dimension of maximal commutative subalgebras of $\gog$: if $d$
denotes the dimension of such a  subalgebra, we have
$$d = \dim\gom_{k(\ov \a)}(\ov\a).$$
where $\ov\a$ is a long simple root having maximal  distance from the
affine node in the extended Dynkin diagram of $\D$.
\smallskip

Then we show  that, when $\a\in\Pi_\ell$,  $\gom_{k(\a)}(\a)$ is a
maximal abelian ideal of $\frak b$ and that each maximal abelian ideal
of $\frak b$ occurs in this way for a unique $\a\in\Pi_\ell$. Hence
the map $\a\mapsto \gom_{k(\a)}(\a)$ is a bijection between $\Pi_\ell$
and maximal ideals in $\II$; the same bijection has been shown via a 
case by case verification in \cite{PR} and in a uniform way in
\cite{P, Corollary 3.8}.
\smallskip

Finally, we find a general dimension formula for all the ideals
$\gom_h(\a)$. In the above notation, to each pair $(\a,h),\,1\leq h<
{k(\a)},\,\a\in\Dp_\ell$  we associate a finite irreducible subsystem
$\D_{h,\a}$ of $\Da$ and we prove that
$$\dim(\gom_h(\a))=g-1+\sum\limits_{i=1}^{h-1}\left(g_h(\a)-1\right)$$
where $g_h(\a)$ is the dual Coxeter number of $\D_{h,\a}$. 
\par

The paper is organized as follows. In section 2 we recall some basic
facts about affine root systems and Coxeter groups. In  section 3 we
recall our encoding of $ad$-nilpotent ideals of $\gog$ via the affine
Weyl group and we prove the two results on $\I$ stated above. In
section 4 we characterize the maximal abelian ideals of $\gob$ in
terms of the corresponding elements in $\widehat W$. In section 5 we
give a structure theorem for the elements of $X$. In section 6 we
introduce $\II(\t)$, we describe its unique minimal element and begin
the study of the maximal elements, proving that they should belong to
$\II(\t,\a)$ with $\a\in\Pi_\ell$. Section 7 is devoted to the proof
of the existence of unique minimal and maximal elements in 
$\II(\t,\a)$. In section 8 we give the explicit description of the
maximal elements $\gom_h(\a)$ which leads to our dimension formulas.
The last section deals with the explicit computation  of
$\dim(\gom_h(\a))$ for any admissible pair $(h,\a),\,\a\in\Pi_\ell$.
\bigskip

\ni{\bf Notational conventions.} 
If $R$ is a root system we shall denote by $\Dp(R),\,Q(R)$
$Q^+(R),\,W(R)$ the positive roots, the root lattice,  the positive
root lattice and the Weyl group of $R$,  respectively. If $Y$ is a
subset of $R$, we shall also  use the notation
$\D(Y)$, and $W(Y)$, to denote the subsystem generated by $Y$, and its
Weyl group, respectively. If $R$ is irreducible, we set $\th(R)$ to be
its highest root. We  reserve the symbols $\Dp,\,Q,\,Q^+,\,\th,\, W$
for the same data relative to the root system $\D$ of $\gog$.\par If
$R$ is simply laced we shall regard all roots as  long.
\par 
The maximal irreducible subsystems of a root system $R$ will be called
the {\it components} of $R$.
\par 
For $\a\in Q,\,\a=\sum_{i=1}^n\lambda_i\a_i$ we set $\supp(\a)= 
\{\a_i\mid\l_i\ne 0\}$.
\par 
We shall mainly use two partial orderings throughout the paper, both
denoted by $\leq$: we shall endow $Q$ (or $\Q$) with the usual order:
$\a<\b$ if $\b-\a$ is a sum of positive roots (resp. coroots). Further
we shall regard $W$ and $\widehat W$ as posets via  the weak order
(more precisely the weak left Bruhat order with respect to the
standard set of Coxeter generators, see \cite{H, 5.9}). Subsets of $W$
or $\widehat W$ will be endowed with the restriction of $\leq$.  Only
in one instance we shall use the Bruhat order: it will be denoted by
$<_B$.
\par 
We shall denote by $A\meno B$ the set difference of $A$ and $B$,
and by $A+B$ the symmetric difference of the sets $A,\, B$. For
$k\in\Bbb N^+$, we set $[k]=\{1,\ldots,k\}$.
\bigskip

\heading 
\S2 Affine root systems and combinatorics of roots
\endheading
\bigskip

Our basic tools are the standard theory of affine Weyl groups and some
combinatorial properties of root systems. This section is devoted to
recollect these topics.
\par
We set $V\equiv \frak h_{\Bbb R}^*=\bigoplus\limits_{i=1}^n \real
\a_i$ and denote by $(\ ,\ )$ the positive  symmetric bilinear form
induced on $V$ by the Killing form. The affine root system associated
to $\D$ can be introduced as follows \cite{Kac, Chapter 6}. We extend
$V$ and its inner product setting $\widehat V= V \oplus \real \d\oplus
\real \l$, $(\d,\d)=(\d,V)=(\l,\l)=(\l,V)=0$, and $(\d,\l)=1$. We
still denote by $(\ ,\ )$ the resulting (non-degenerate) bilinear
form. The affine root system associated to $\D$ is $\Da= \D+\ganz
\d=\{\a+k\d\mid \a\in \D, \ k\in \ganz\}$; remark that the affine
roots are non isotropic vectors. The set of positive affine roots is
$\Dap=(\Dp+\nat \d)\cup (\Dm+\nat^+\d)$, where $\Dm=-\Dp$. We denote
by $\theta$ the highest root of $\D$ and set $\a_0=-\theta+\d$,
$\widehat \Pi=\{\a_0,\dots,\a_n\}$. We sometimes call the roots of
$\D$ inside $\Da$ the \lq\lq finite roots\rq\rq. For each $\a\in\Dap$
we denote by $s_\a$ the corresponding reflection of $\widehat V$,
$s_\a(x)=x-{{2(\a,x)}\over{(\a,\a)}} x$ for  $x\in \widehat V$. The
affine Weyl group associated to $\D$ is the group $\widehat W$
generated by $\{s_\a\mid \a\in \Dap\}$; the set $\{s_\a\mid \a\in
\widehat\Pi\}$ is a set of Coxeter generators for $\widehat W$;  we
denote by $\ell$ the corresponding length function. Note that
$w(\d)=\d$ for any $w\in\widehat W$. It turns out that $\widehat W$ is
a semidirect product $T\rtimes W$, where $T=\{t_\tau \mid \tau\in\Q\}
\cong \Q$ is the subgroup of {\it translations}, and the action of $W$
on $T$ is $vt_\tau v^{-1}=t_{v(\tau)}$. The action of $t_\tau$ on
$V\oplus \real \delta$, in particular on the roots, is given by
$t_\tau(x)=x-(x, \tau) \delta$, for each $x\in V\oplus \real \delta$.
(See \cite{Kac} for the general definition of $t_\tau$ on $\widehat
V$).
\par
Now consider the $\widehat W$-invariant affine subspace $E=\{x\in
\widehat V\mid (x,\d)=1\}$. Remark that $E=V\oplus\Bbb R\d+\l$; let $\pi:E\to
V$ be the natural projection. For $w\in \widehat W$ we set $\overline
w=\pi\circ w_{|E}$. We have that   the map $w\mapsto \overline w$
gives an isomorphism of $\widehat W$ onto a group  $\Waf$ of affine
transformations of $V$ \cite{Kac, 6.6}. Under this isomorphism  each
element $t_\tau$ of $T$ maps to the true translation by $\tau$. So
$\Waf$ is the usual affine representation of the affine Weyl group
\cite{B, VI, \S 2}.
In fact, for $\a\in \Dp$ and $k\in \ganz$ set $$H_{\a,k}=\{x\in V\mid
(x,\a)=k\}$$ and let $s_{\a,k}$ denote the orthogonal reflection with
respect to $H_{\a,k}$ in $V$. For $\b\in\widehat V$, let  $\widehat
H_\b$  be the hyperplane orthogonal to $\b$ in $\widehat V$. Then it
is easy to see that, for each $\a\in \Dp$, $k> 0$ and $h\geq 0$, we
have
$$\align 
\pi(\widehat H_{-\a+k\d}\cap E)= H_{\a,k}, \quad\quad  
&\overline{s_{-\a+k\d}}=s_{\a,k},\\ 
\pi(\widehat H_{\a+h\d}\cap E)= H_{\a,-h}, \quad\quad  
&\overline{s_{\a+h\d}}=s_{\a,-h}.\\ 
\endalign$$
Hence we have a transparent correspondence between the linear and the 
affine context.  
Henceforward we shall omit the bar and use the same notation for 
the elements in $\hatw$ and the corresponding elements $\Waf$.
We denote by $C$  is the fundamental chamber of $\widehat W$  and
by $C_1$  the fundamental alcove of $\Waf$, 
$$\align 
C&=\{x\in \widehat V\mid (x,\a_i)>0 \text{ for } i=0,\dots ,n\},\\
C_1&=\{x\in V\mid (x,\a)>0\,\forall \a\in\Pi,\ (x,\theta)<1\}\\
&=\{x\in V\mid 0<(x,\a)<1\,\forall \a\in\Dp\}.
\endalign$$ 
Then we have $C_1=\pi(C\cap E)$. This implies the following result.
\bigskip

\proclaim{2.1} 
For any $\a\in\Dp,\,k\in\Bbb N^+,\,h\in\Bbb N$ $$ {w^{-1}(-\a+k\d)<
0\text{ if and only if $H_{\a,k}$ separates }  C_1 \text { and }
w(C_1),}\atop {w^{-1}(\a+h\d)< 0 \text{ if and only if } H_{\a,-h}
\text{ separates } C_1 \text{ and } w(C_1).}
$$\endproclaim
\bigskip

\ni
{\bf 2.2.\ } 
For $j>0$ set $$C_j=\{x\in V\mid (x,\a)>0\,\forall \a\in\Pi,\
(x,\theta)<j\};$$ also set $$C_\infty=\{x\in V\mid (x,\a_i)>0\text{
for }i=1,\ldots,n\}.$$ Thus $C_\infty$ is the fundamental chamber of
$W$. The walls of $C_j$ are the hyperplanes, $H_{\a_i,0},\,1\leq i\leq
n$ and $H_{\theta,j}$,  so that $C_j=jC_1$. For all $\b\in \Dp$ and
$x\in C_\infty$ we have $(x,\b)\leq (x,\th)$. In particular we obtain
that if $H_{\b,k}\cap C_2\ne\emptyset$, then $k=1$.
\bigskip

\ni
{\bf 2.3.\ } 
Let $R$ be a finite or affine  root system. Choose the set  $R^+$  of
positive roots, and denote by $R_\Sigma$ the corresponding set of
simple roots. Let $G$ be the Weyl group of $R$; for $\a\in R^+$ we
denote by $s_\a$ the reflection with respect to $\a$.  For $v\in G$ we
set $$N(v)=\{\b\in R^+\mid v^{-1}(\b)<0\}.$$ Moreover, we denote  by 
$L_v$ the set of left  descents of $v$: 
$$ L_v=\{\a\in R_\Sigma\mid\ell(s_\a v)<\ell(v)\}.$$ 
Similarly we define the set of {\it right descents} of $v$ as $R_v=\{\a\in R_\Sigma\mid\ell(vs_\a) <\ell(v)\}.$\par
The following facts are well known; proofs and details about
(3) and (4) can be found in  \cite{D}, \cite{CP1}. Item (5) follows from (3).
\item{(1)}   
$L_v=N(v)\cap R_\Sigma$;
\item{(2)} 
$\a\in L_v$ if and only if there exists a reduced expression of $v$
starting with $s_\a$.
\item{(3)} 
$N(v)$ is finite, closed under root addition and, moreover, if $\b_1, \dots,
\b_k\in R^+$, $c_1, \dots, c_k\in \Bbb R^+$, and
$c_1\b_1+\cdots+c_k\b_k\in N(v)$, then $\b_i\in N(v)$ for at least one
$i\in \{1, \dots, k\}$.
\item{(4)} 
Conversely, if $R\not\cong\widehat A_1$, and  $N\sbsq R^+$ is a
finite set which satisfies the following conditions:
\roster
\item"({\it a})" 
$N$ is closed under root addition, 
\item"({\it b})"
for all $\b_1, \b_2\in R^+$  such that $\b_1+\b_2\in N$,  
at least one of $\b_1, \b_2$ belongs to\nobreak\ $N$,
\endroster
\item{}
then there exists $v\in W(R)$ such that $N=N(v)$. 
If $R\cong\widehat A_1$, then the sets $N(v)$, $v\in W(R)$, are exactly 
all the sets $\{-\a+i\d\mid 1\leq i\leq k\}$, for all $k\in \Bbb N^+$, and 
$\{\a+i\d\mid 0\leq i\leq k\}$, for all $k\in \Bbb N$, where $\a$ is the 
finite simple root.\smallskip
\item{(5)} Let $I\subseteq R_\Sigma$. Then $L_v\subseteq I$ if and
only if for all $\a\in N(v)$ we have $\supp(\a)\cap I\ne\emptyset$. 

\bigskip

We describe explicitly the set $N(v)$. Part (1) of the following statement is
clear; part (2) is  well-known \cite{B, VI, 1.6}.
\bigskip
 
\proclaim{2.4}
Let $w\in G$ and let $w=\sigma_1\cdots \sigma_k$ be any reduced
expression of $w$ with $\sigma_i=s_{\b_i}$, $\b_i\in R_\Sigma$. If we
set $\g_1=\b_1,\g_2=\sigma_1(\b_2),\ldots, \g_k=\sigma_1\cdots
\sigma_{k-1}(\b_k)$, then we have:
\roster
\item
$w=s_{\g_k}\cdots s_{\g_2} s_{\g_1}$.
\item
$N(w)=\{\g_1,\ldots,\g_k\}$.
\endroster
\endproclaim
\bigskip

\ni{\bf 2.5}
We also set $\Npm(v)=N(v)\cup-N(v)$. It is easy to see  that for any
$u,\, v\in G$ we have $$\Npm(vu)=\Npm(v)+v\Npm(u),$$ where we denote
by $+$ the symmetric difference of sets. 
\bigskip

We recall that the weak order $\weak$ of $G$ is defined as follows:
for $v, v'\in G$, we have $v\leq v'$ if and only if there exists an
element $u\in G$ such that $v'=vu,\,\ell(v')=\ell(v)+\ell(u)$. From
2.4 and 2.5 we easily obtain the following characterization of the 
weak order.
\smallskip
 
\proclaim{2.6}
$v\weak v'\text{ if and only if } N(v)\sbsq N(v').$
\endproclaim
\bigskip
Moreover, we have the following results.
\smallskip
\proclaim{2.7} 
The following conditions $(1)$, $(2)$,  $(3)$ and $(4)$ are equivalent: 
\parni
$(1)$ $N(vu)\spsq N(v)$; $(2)$ $v N(u)\sbsq \Dp$; $(3)$
$\ell(vu)=\ell(v)+\ell(u)$; $(4)$ $v\leq vu$. 
\par
Similarly, the following conditions $(1')$, $(2')$,  $(3')$ and $(4')$ are
equivalent:
\parni
$(1')$ $N(vu)\sbsq N(v)$; $(2')$ $v N(u)\sbsq \Dm$; $(3')$ $\ell(v
u)=\ell(v)-\ell(u)$; $(4')$ $v\geq vu$.
\endproclaim
\bigskip

Let $G'$ be a standard parabolic subgroup of $G$, and let $R'$ be its
root system. Recall that any $v\in G$ can be decomposed as
$v=v'v''$ where $v'\in G'$ and $v''$ is an element of minimal length
in the right coset $G'v$; moreover, $v',\,v''$ are unique. It follows that
\smallskip

\proclaim{2.8}
For all $v\in G$ there exists  $v'\in G'$ such that $N(v)\cap
R'=N(v')$.
\endproclaim
\smallskip

Set moreover $R''=R\meno R'$. It is well known that  $y\in G'x$ is the
minimal length representative of the right coset $G'x$ if and only if
$L_y\sbsq R''$. Therefore, $\{y\in G\mid L_y\sbsq R''\}$ is the set of
minimal length representatives for the set of right cosets
$G'\backslash G$. Thus the restriction of the weak order of $G$ on 
$\{y\in G\mid L_y\sbsq R''\}$ induces a partial ordering on
$G'\backslash G$. When saying {\it the poset  $G'\backslash G$}, we
shall always refer to this ordering.
\bigskip

\ni{\bf 2.9.}
The poset $G'\backslash G$ has a unique minimal and a unique maximal
element. The identity $1$ of $G$  clearly corresponds to the minimum of
$G'\backslash G$. If $w_0$ is the longest element of of $G$ and $w_0'$
is the longest element of $G'$, then we have that $N(w_0'w_0)=
\Dp(R)\meno \Dp(R')$. By 2.3 (5), if $y\in G$ is such that
$L_y\sbsq R''$, then $N(y)\sbsq \Dp(R)\meno \Dp(R')$; therefore we
obtain that $w_0'w_0$ is the unique maximal element in $\{y\in G\mid
L_y\sbsq R''\}$, with respect  to the weak order. Thus $w_0'w_0$
corresponds to the unique maximal element of  $G'\backslash G$.

\vskip 1cm

\heading \S3 
$ad$-nilpotent ideals via affine Weyl groups
\endheading
\bigskip

We turn now to $ad$-nilpotent ideals. Recall from \cite{CP2, \S2} that
to any $ad$-nilpotent ideal $\goi=\sum\limits_{\a\in\p}\gog_\a$ we can
associate  a set of positive roots in the affine root system attached
to $\gog$ in the following way: $$L(\goi)=\bigcup_{k\geq 1}\left(
-\p_\goi^k+k\d\right),$$ where $\p_\goi^1=\p$ and
$\p_\goi^k=(\p_\goi^{k-1}+\p)\cap \D,\,k\geq 2$. In \cite{CP2, Prop. 2.4}
we proved that there exists a (unique) $w_\goi\in \hatw$ such that  
$N(\wi)=L(\goi)$. Set $$\Cal
W=\left\{w_\goi\mid\goi\in\I\right\},\qquad\Cal
W_{ab}=\left\{w_\goi\mid\goi\in\I_{ab}\right\},$$ so that
$\goi\mapsto w_\goi$ is a bijection $\I\to\Cal W$, restricting to a
bijection $\II\to\Cal W_{ab}.$ The inverse map will be denoted by
$w\mapsto \goi(w)$.\par The sets $\Cal W,\,\Cal W_{ab}$ have been
characterized in our previous papers.
\bigskip

\proclaim{Proposition 3.1} 
$(1)$ \cite{CP3, Prop. 2} Let $w\in \widehat  W$, $w=t_\tau v$,
$\tau\in \Q$, $v\in W$. Then $w\in\Cal W$  if and only if the
following conditions hold:
\roster
\item"$(i)$" $w(C_1)\subset C_\infty$; 
\item"$(ii)$" $(v^{-1}(\tau),\a_i) \leq 1$ for each $i\in [n]$ and
$(\tau, v(\theta)) \geq - 2$.
\endroster 
\ni
$(2)$ \cite{CP2, Th. 2.9} $w\in\Cal W_{ab}$ if and
only if $w(C_1)\subset C_2$.
\endproclaim
\bigskip

An obvious but useful observation is the following remark (see 2.6).
\bigskip

\proclaim{3.2} 
For $\goi,\goj\in \I$, the following conditions are equivalent: \parni
$$(1)\ \goi\subseteq \goj;\qquad (2)\  N(w_\goi)\subseteq N(w_\goj);
\qquad (3)\  w_\goi\leq w_\goj.$$
\endproclaim 
\smallskip

Finally, we remark that the construction given at the beginning of the
section is quite simple in the special case of abelian ideals. Indeed,
we have the following result.

\bigskip
\proclaim{Proposition 3.3}
{\rm (Peterson, see also \cite{CP2, section 2})}
Let $\goi\in \II$. Then $N(w_\goi)=\{-\a+\d\mid \gog_\a \sbsq \goi\}.$
In particular, $dim(\goi)=\ell(w_\goi)$.
\par 
Conversely, assume that $w\in \widehat W$ is such that $N(w)\sbsq
\Dm+\d$ and set $\goi=\sum\limits_{-\a+\d\in N(w)}\gog_\a$. Then
$\goi\in  \II$.
\endproclaim
\bigskip 

The  next results complete the analysis developed in \cite{CP2},
\cite{CP3}, showing once again that all the combinatorial information
regarding $ad$-nilpotent ideals  can be recovered from $w_\goi$.
\parni
Remark that $\sum\limits_{\a\in\p}\gog_\a\in\I$ if and only if
$\p$ is a dual order ideal of $\Dp$ endowed with the usual partial
order, i.e., if and only if, for all $\a\in \Phi$ and $\b\in \Dp$
such that $\a\leq \b$, we have that $\b\in \Phi$. 
\par
For all $\goi\in \I$ we set 
$$\p_\goi=\{\b\in\Dp\mid \gog_\b\sbsq\goi\},$$
and denote by $\Cal A_\goi$ the set of minimal elements of $\p_\goi$.
It is clear that $\Cal A_\goi$ is an antichain of $\Dp$, i.e. a set of
incomparable elements. As explained in \cite{CP3, \S4 (1)}, $\p_\goi$
is completely determined by $\Cal A_\goi$.
\bigskip

\proclaim{Proposition 3.4}
For $\goi\in\I$ the following statements are equivalent:
\roster
\item 
$\wi(\a_i)=\b-\d,\,\b\in\Dp$ for some $i,\,0\leq i\leq n$;
\item 
$\b\in \Cal A_\goi$.
\endroster
\endproclaim

\demo{Proof}
$(1)\implies (2)$. We have $\wi(\a_i)=\b-\d<0$; hence, 
$\wi s_i(\a_i)=-\b+\d>0$. By 2.7, it follows that $N(\wi)
=N(\wi s_i)\cup \{-\b+\d\}$. If there exists $\g\in \Phi_\goi$ 
such that $\b>\g$, then $-\g+\d$ is equal to $-\b+\d$ plus a 
sum of roots in $\Dp$. But $-\g+\d\in N(\wi)$, hence 
$-\g+\d\in N(\wi s_i)$, and since $-\b+\d\not\in N(\wi s_i)$, 
by 2.3 (3) we obtain that some root in $\Dp$ belongs to $N(\wi s_i)$. 
This is not possible, since $N(\wi s_i)\sbsq N(\wi)\sbsq \Dm+\nat^+\d$.
\parni
$(2)\implies (1)$. 
Let $\b\in \Phi_\goi$ be minimal. We shall prove
that $L(\goi)\meno\{-\b+\d\}$ verifies conditions ({\it a}) and 
({\it b}) of 2.3 (4). This will imply that there exists $w\in \hatw$ such that 
$L(\goi)\meno\{-\b+\d\}=N(w)$; hence, by 2.6 and 2.5, $\wi=w s_i$ 
for some $i\in \{0,\dots,n\}$ and  $\wi(\a_i)=-w(\a_i)=\b-\d$.
Set $N=L(\goi)\meno\{-\b+\d\}$. 
Since $N\sbsq \Dm+\nat^+\d$, no pair of roots in $N$ sum up
to $-\b+\d$; since, moreover, $L(\goi)$ is closed under root addition we
obtain that $N$ is closed, too.  
It remains to prove that ({\it b}) holds.  Assume $\eta\in N$ and
$\eta=\eta_1+\eta_2$, with $\eta_1, \eta_2\in\Dap$. If $\eta_1, \eta_2\ne
-\b+\d$, we are done since ({\it b}) holds for $L(\goi)$.
We assume that $\eta_1=-\b+\d$ and prove that then $\eta_2\in N$.  
Set $\eta=-\ov\eta+k\d$, $k>0$, $\ov\eta\in \Dp$. 
We prove the claim by induction on $k$. 
First notice that $\eta_2\not\in\Dp+\nat \d$, since, 
if $\eta_2=\a+(k-1)\d$ with $\a\in\Dp$, 
then $-\ov\eta=-\b+\a$, which implies $\b>\ov\eta$ against the
minimality of $\b$ in $\p_\goi$.
So we may assume $\eta_2=-\a+(k-1)\d$, $k>1$, $\a\in\Dp$.
We first assume $k=2$, so that $\eta=(-\b+\d)+(-\a+\d)$. 
By assumption $\ov\eta\in\p_\goi^2$, hence $\ov\eta=\b'+\a'$ 
with $\b', \a'\in \p_\goi$. Moreover we have that $\ov\eta=\b+\a$, 
hence $0<(\a'+\b', \a+\b)=(\a', \a)+(\a', \b)+(\b', \a)+(\b', \b)$. 
It follows that at least one of the scalar products in the right hand
side of the above identity is positive. If $(\a', \a)>0$, then
$\a'-\a\in \D$, and since $\b=\b'+(\a'-\a)$, we obtain that
$\a'-\a<0$, by the minimality of $\b$ in $\p_\goi$. Since
$-\a'+\d=-\a+\d+(\a-\a')$, $-\a'+\d\in L(\goi)$, $\a-\a'\in
\Dp\meno L(\goi)$, and $L(\goi)$ has property ({\it b}), we obtain
that $-\a+\d\in L(\goi)$, hence $-\a+\d\in N$.
If $(\b', \a)>0$, we obtain similarly that $-\a+\d\in N$. 
So we assume that $(\a', \b)>0$. Then by minimality of $\b$ we have
that $\a'-\b\in \Dp$, hence $\a-\b'=\a'-\b\in\Dp$ and we can argue 
as in the previous case, so obtaining that  $-\a+\d\in N$. The case   
$(\b', \b)>0$ is similar. Thus we have proved the claim for $k=2$. 
Now we assume $k>2$ and the claim true for $k-1$. 
By assumption we have $\ov\eta=\b'+\a'$ with $\a'\in \p_\goi^{k-1}$ 
and $\b'\in \p_\goi$. Arguing as above, we find that at least one of 
$(\a', \a)$, $(\a', \b)$, $(\b', \a)$, $(\b',\b)$, is positive.
If $(\a',\a)>0$, then by minimality of $\b$ we obtain that 
$\a-\a'>0$; moreover $-\a'+(k-1)\d=-\a+(k-1)\d+(\a-\a')$, and 
since $-\a'+(k-1)\d\in L(\goi)$ and $\a-\a'\notin L(\goi)$ we obtain 
that $-\a+(k-1)\d\in L(\goi)$, hence $-\a+(k-1)\d\in N$. 
If $(\b', \b)>0$, we obtain as above that $\a-\a'>0\in \Dp$, and 
we are reduced to the previous case. So we assume that 
$(\a', \b)>0$. As above, $\a'-\b>0$, hence  there exists $\g\in \Dp$ such that $\a'=\b+\g$
and $\a=\b'+\g$. Hence we have that $-\a'+(k-1)\d=(-\b+\d)+(-\g+(k-2)\d)$. 
By the induction assumption we obtain that $-\g+(k-2)\d\in L(\goi)$, 
hence $-\a+(k-1)\d=(-\b'+\d)+(-\g+(k-2)\d)\in L(\goi)$ and therefore 
$-\a+(k-1)\d\in N$. The case $(\b', \a)>0$ is similar.  
\endemo
\bigskip

\ni 
For each $\b\in \Cal A_\goi$, let $x_\b$ be any  non zero element in
$\gog_\b$. Then $\{x_\b\mid \b\in \Cal A_\goi\}$  is a minimal set of
generators of $\goi$ as a $\gob$-module. Hence Proposition 3.4,
together with \cite{CP2, Prop. 2.12, ii}, imply the following result.
\bigskip

\proclaim{Corollary 3.5}
For any $\goi\in\I$, the minimal number of generators of $\goi$ as a
$\gob$-module equals the number of right descents of $w_\goi$ and the
cardinality of the antichain $\Cal A_\goi$.
\endproclaim
\bigskip

To proceed further, we need to recall one of the main results from
\cite{CP3}. Set 
$$D=\left\{\sigma\in  \Q \mid (\sigma,\alpha_i)\leq 1 \text{ for each
} i\in \{1,\dots ,n\} \text{ and } (\sigma,\theta)\geq -2\right\}.$$
Set moreover 
$$\dab=\{\sigma\in  \Q \mid (\sigma,\alpha)\in\{0,1,-1,-2\} \text{ for each
} \a\in\Dp\}.$$
It is clear that $\dab\sbsq D$. 
The first part of next statement is a result of \cite{CP3},
generalizing  the second one, which appears in \cite{K2}. \bigskip

\proclaim 
{Proposition 3.6}
\cite{CP3, Prop. 3}, \cite{K2, Theorem 2.6}  
The map $F: t_{\tau} v\mapsto v^{-1}(\tau)$ is a bijection between 
$\Cal W$ and $D$.
Moreover, $F$ restricts to a bijection between $\wab$ and $\dab$.
\par
\endproclaim 

\nd 
Set now, for $1\leq j\leq ht(\theta)$: $$D'_j=\left\{\sigma\in \Q \mid
-(j+1)\leq (\sigma,\b)<j+1\,\forall\,\b\in\Dp\right\}.$$ Set
$D_i=D'_i\cap D$. Note that  $D_{ab}=D_1=D_1'$. Next lemma follows by
standard arguments (see \cite{IM}); for the sake of completeness we
include a proof.
\bigskip

\proclaim{Lemma 3.7} 
Assume  $\goi\in\I$ and $w_{\frak i}=t_\tau v$, with $\tau\in \Q$ and
$v\in W$. The following conditions are equivalent:
\roster
\item"$(i)$" 
$w_{\frak i}(C_1)\subset C_{j+1}$;
\item"$(ii)$" 
for each $\alpha \in \Delta^+$,   $0\leq (\tau, \alpha)\leq j+1$;
moreover, if $(\tau, \alpha)=0$, then $v^{-1}(\alpha)>0$ and if
$(\tau, \alpha)= j+1$, then $v^{-1}(\alpha)<0$.
\endroster
\endproclaim

\demo{Proof}
Assume that $(i)$ holds; then   the hyperplanes
$H_{\a,0},\,H_{\a,j+1}$ do not separate $C_1$ and $w_\goi(C_1)$ for
any $\a\in\Dp$. By direct computation, we obtain
$$w_\goi^{-1}(-\a+(j+1)\d)=v^{-1}t_{-\tau}(-\a+(j+1)\d)=-v^{-1}(\a)+(j
+1-(\a,\tau))\d.$$ 
By 2.1, the above root should be positive, hence  $(\a,\tau)\leq j+1$,
and, if equality holds, then necessarily $v^{-1}(\a)<0$. Similarly, 
$w_\goi^{-1}(\a)=v^{-1}(\a)-(\a,\tau)\d>0$, hence $(\a,\tau)\geq 0$,
and $v^{-1}(\a)>0$ if $(\a,\tau)=0$. Hence $(ii)$ holds. The converse
is similar.
\endemo
\bigskip

For $\goi\in \I$, let $n(\goi)$ denote the nilpotence index of
$\goi\in\I$, i.e. the minimal $n\in \nat^+$ such that $\goi^{n+1}=0$, where
$\goi^n$ is the $n$th term of the descending central series of
$\goi$ ($\goi^1=\goi$).

\proclaim{Proposition 3.8} 
$F$ induces by restriction  bijections between $\{w_\frak i \mid \frak
i\in \I,\, n(\goi)\leq j \}$ and $D_j$.
\endproclaim

\demo{Proof} By  \cite{CP2, Corollary 2.5}  the ad-nilpotent ideal
$\frak i$ of $\frak b$ has nilpotence index at most $j$ if and only
if $-\a+(j+1)\d\notin L_\goi$ for any $\a\in\Dp$; in turn, this
condition is equivalent to $$w_{\frak i}(C_1)\subset C_{j+1}$$ (this
follows, e.g., by formulas 2.1). By Lemma 3.7,  we assume that  
condition $3.7 (ii)$  holds and consider $(v^{-1}(\tau),\beta)$ for
$\beta\in \Delta^+$.  We have either $\beta=v^{-1}(\alpha)$ or
$\beta=-v^{-1}(\alpha)$ for some $\alpha\in \Delta^+$. In the former
case $v^{-1}(\alpha)>0$ and by our assumption $0\leq
(v^{-1}(\tau),\beta)= (\tau, \alpha)<j+1$; in the latter case
$v^{-1}(\alpha)<0$, whence $-j-1\leq (v^{-1}(\tau),\beta)= -(\tau,
\alpha)<0$. In either cases we have $-j-1\leq (v^{-1}(\tau),\beta)
<j+1$. In a similar way we see that, conversely,  if $-j-1\leq
(v^{-1}(\tau),\beta) <j+1$, then $3.7 (ii)$ holds.
\endemo
\vskip 1cm

\heading \S4 
A characterization of maximal abelian ideals
\endheading
\bigskip

In this section we study maximal abelian ideals of $\gob$, obtaining a
characterization for them in terms of $\wi$.
\par

We shall need the following easy remark.
\bigskip

\ni
{\bf 4.1.} If $\t\in \Q$ is dominant, then either $(\t,\th)=0$, or
$(\t,\th)\geq 2$. This is immediate since $0$ is  the unique element
in $\Q\cap\overline C_1$.
\bigskip

\proclaim{Lemma 4.2}
Let $\goi\in\II$ and $\b\in \D$. If $\b+k\d\in w_\goi(\widehat \Pi)$,
for some $k \geq 2$, then $\b=-\th$, $k=2$, and $w_\goi^{-1} (-\th+2\d)
\in \Pi$.
\endproclaim

\demo{Proof}
Let $\b+k\d= w_\i(\a)$,  $\a\in \widehat\Pi$. Then $\a\not\in
L_{w_\goi^{-1}}$, hence $\ell(w_\goi s_\a)>\ell(w_\goi)$ and by 2.7
$N(w_\goi s_\a)=N(w_\goi)\cup\{w_\goi(\a)\}=N(w_\goi)\cup\{\b+k\d\}$.
If $\b>0$, then by 2.3~(3) we should obtain that $\b\in N(w_\goi s_\a)$,
which is impossible since $\b\notin N(w_\goi)$; hence $\b<0$. In 
particular,  $\b+k\d\notin \Pi+\nat \d$, hence, by 3.3, no root in $\Pi+\nat \d$ 
belongs to $N(w_\goi s_\a)$.
But it is clear that $\b+k\d$ is equal to $-\th+2\d$ plus a linear
combination of roots in $\Pi+\nat \d$ with integral non-negative
coefficients; hence by 2.3 (3) we obtain that $-\th+2\d\in N(w_\goi s_\a)$.
Using 3.3 again, we obtain that $-\th+2\d=\b+k\d$. 
It remains to prove that $\a\not=-\th+\d$. Let $w_\goi=t_\t v$, $\t\in
\Q$, $v\in W$. If by contradiction $\a=-\th+\d$, then $v(\th)=\th$ and
$t_\t v(-\th+\d)= -\th+(1+(\th, \tau)) \d$. Since $\t$ is dominant, by
4.1 we obtain that either $t_\t v(-\th+\d)=-\th+\d$ or $t_\t
v(-\th+\d)=-\th+k\d$ with $k>2$.
\endemo
\bigskip

\proclaim{Lemma 4.3}
Let $\goi\in \II$. Then $w_\goi(\widehat \Pi)\sbsq \pm(\Dm+\d)\cup
\Pi\cup\{-\th+2\d\}$.
\endproclaim

\demo{Proof}
Let $\a\in \widehat \Pi$. If $w_\goi(\a)<0$, then by 3.3 $w_\goi(\a)
\in \Dp-\d$. Assume $w_\goi(\a)>0$, so that, by 2.7, $N(w_\goi
s_\a)=N(w_\goi)\cup\{w_\goi(\a)\}$.
If $w_\goi(\a)=\b+k\d$ with $\b\in \Dp$ and $k\geq 0$, then by 2.3 (3)
there exists some $\g\in \Pi$ such that $\g\in N(w_\goi s_\a)$,
whence, by 3.3, $\g=w_\goi(\a)$. Finally, if $w_\goi(\a)=-\b+k\d$ with
$\b\in \Dp$ and $k>0$, then either $k=1$, or, by 4.2, $\b=\th$ and
$k=2$.
\endemo
\bigskip

\proclaim{Theorem 4.4} 
$\goi\in\II$ is maximal if and only if $w_\goi(\Pia)\cap
(\Dm+\d)=\emptyset$. In this case $-\th+2\d\in w_\goi(\Pi_\ell)$.
\endproclaim

\demo{Proof} 
The abelian ideal $\goi$ is maximal in $\wab$ if and only if, for all 
$w\in \hatw$ such that $w_\goi \leq w$ in the weak order, $w\not\in \wab$.
By 3.3, if $u\in \wab$ and $u'\leq u$ in the weak order of $\hatw$, 
then $u'\in \wab$ as well. This implies that $\goi$ is maximal in $\II$
if and only if, for all $\a\in \Pia$ such that $w_\goi(\a)>0$, we have
that $w_\goi s_\a\not\in \wab$. Since for $w_\goi(\a)>0$ we have 
$N(w_\goi s_\a)=N(w_\goi)\cup\{ w_\goi(\a)\}$, this happens if and 
only if $w_\goi(\a)\not\in \Dm+\d$. This proves the first statement.
\par
Now if $w_\goi(\Pia)\cap(\Dm+\d)=\emptyset$, we have in
particular that $w_\goi$ is not the identity of $\widehat W$. Notice,
moreover, that by 4.1 no translation may correspond to some non zero
abelian ideal.  Hence if $w_\goi=t_\t v$, then $v$ is not the
identity; therefore $v(\a)<0$ for at least one $\a\in \Pi$. 
But by 4.3, $t_\t v(\Pia)\sbsq\Dp-\d\cup\Pi\cup\{-\th+2\d\}$, 
hence for such an $\a$, $v(\a)=-\th$ and $\wi(\a)=-\th+2\d$.
Moreover, since $\th$ is long, $\a$ is long too.
\endemo
\bigskip

\ni
{\bf 4.5. Remark.} 
If $w\in\wab$, and $w=t_\t v$ with $\t\in \Q,\,v\in W$, then
$w^{-1}(-\th+2\d)\in \Pi_\ell$ if and only if $v^{-1}(-\th)\in
\Pi_\ell$. In fact, if $v(\a)=-\th$, then by 4.1 we obtain that either
$w(\a)=-\th$, or $w(\a)=-\th+2\d$, but  by 4.3 the first case cannot
occur.
\bigskip

In 3.6 we recalled Kostant's bijection between $\II$ and a certain 
subset $\dab$ of $\Q$. Theorem 4.4 translates into the following 
characterization of maximal abelian ideals inside $\dab$.
\bigskip

\proclaim{Proposition 4.6} 
Let $\goi\in \II$ and let $w_\goi=t_\t v$, with $\t\in \Q$ and $v\in
W$. Then $\goi$ is maximal if and only if
\roster
\item 
$(v^{-1}(\tau),\a_i)\ne -1$ for $1\leq i\leq n$;
\item 
$(v^{-1}(\tau),\th)\ne 0$.
\endroster
\endproclaim

\demo{Proof}
For all $\a\in \widehat \D$ we have  $w_\goi(\a)=v(\a)-(v(\a),\t)\d$;
moreover, $v(-\th+\d)=-v(\th)+\d$ and $(\d,\t)=0$. Therefore by a
direct computation we obtain that $w_\goi(\widehat{\Pi})\cap
\D^-+\d=\emptyset$ if and only if (1) and (2) hold.
\endemo
\vskip 1cm

\heading 
\S5 The  structure of $X$
\endheading
\medskip

As  explained in the introduction,  the set $X=\Q\cap \overline C_2$
plays a prominent role in the description of the abelian ideals. In
this section we provide a technical characterization of $X$.
Indeed we shall give an orthogonal decomposition of each element of
$X$ as a sum of  highest roots of nested irreducible parabolic
subsystems of $\D$.
\smallskip 
We need preliminarily the analysis of which roots are dominant; this
fact is known to the experts, nevertheless we prefer to provide a
proof (which is different from the one given in \cite{St1, Proposition
2.1}).
\par 
It is
clear that the highest root $\th$ is dominant. Moreover,  $\th^\vee$
is a short root in the dual system $\D^\vee$ and is clearly dominant
relatively to $\D^\vee$. Symmetrically, if $\th(\D^\vee)$ is the
highest root of $\D^\vee$ and we set $\th_s={2\th(\D^\vee)\over
{(\th(\D^\vee),\th(\D^\vee))}}$, we obtain that $\th_s$ is a short
dominant root in $\D$. Therefore $\th,\,\th_s$ are  dominant; the
converse is also true.
\bigskip

\proclaim{Lemma 5.1}
The unique dominant roots in $\D$ are $\th$ and $\th_s$.
\endproclaim

\demo{Proof}
Let $\a\in\D$ be a long dominant root. Then $(\a, \b)\geq 0$ for all 
$\b\in \Dp$; indeed, since $\a$ is long, we have $(\b, \a^\vee)\in\{0, 
1\}$ for all $\b\in \Dp$. It follows that $\b-\a\not\in \Dp$ for all 
$\b\in \Dp$, otherwise $(\b-\a,\a^\vee) =(\b,\a^\vee)-(\a, 
\a^\vee)<0$, against the assumption that $\a$ is dominant. Therefore 
$\a=\th$. Now assume that $\a$ is a short dominant root. Then 
$\a^\vee$ is a long dominant root in $\D^\vee$, therefore 
$\a^\vee=\th(\D^\vee)$ and $\a=\th_s$.
\endemo
\bigskip

\proclaim{Lemma 5.2}
If $\t\in Q$ is dominant and non zero then $\t\geq \th_s$.
Equivalently, if $\t\in Q^\vee$ is dominant and non zero then
$\t\geq \th^\vee$.
\endproclaim

\demo{Proof}
It is well known that the fundamental weights are nonnegative linear 
combinations of $\Pi$, hence any dominant weight $\s$ is a 
nonnegative linear combination of $\Pi$; it is also clear that 
$\supp(\s)=\Pi$. In particular, if $\s\in-Q^+$, then there 
exists some $\a\in\Pi$ such that $(\s, \a)<0$. Indeed, 
if $0\ne\s=\sum\limits_{i=1}^n \l_i\a_i$ with $\l_i\leq 0$, 
then there exists an $i$ such that $\l_i<0$ and $(\s, \a_i)<0$. 
This is immediate, since $(\a_i,\a_j)\leq 0$ for $1\leq i<j\leq n$, 
hence $(\s, \a_j)\geq 0$ for all $j$ such that $\l_j=0$. 
\par
Now let $\t\in Q$ be dominant and assume, by contradiction, that 
$\t-\th_s\not\in Q^+$. Then $\t-\th_s$ is an integral linear 
combination of $\Pi$ with some negative coefficient, say 
$\t=\th_s-\sum_{j=1}^h \l_j\a_{i_j}+ \sum_{j=h+1}^n\l_j 
\a_{i_j}$, with $\{i_1, \dots ,i_n\} =[n]$, $1\leq h\leq n$, $\l_j\in 
\nat$ and $\l_j>0$ for some $1\leq j\leq h$.
Set $\s=\sum_{j=1}^h \l_j\a_{i_j}$, and $\u=\sum_{j=h+1}^n\l_j 
\a_{i_j}$
Let $\s'\in Q^+$ be a maximal element, with respect to the 
standard partial order, such that $\s'\leq\s$, 
and $\th_s-\s'\in \D$. Set $\th_s'=\th_s-\s'$ and 
$\s''=\s-\s'$, so that 
$\t=\th_s'-\s''+\u$. By the above discussion, if $\s''\ne 0$ 
then there exist some $k\in [h]$ such that $\s''>\a_{i_k}$ and 
$(-\s'', \a_{i_{k}})<0$. 
But it is clear that $(\u,\a_{i_{k}})\leq 0$; 
moreover, $(\th_s', \a_{i_{k}})\leq 0$, otherwise 
$\th_s'-\a_{i_{k}}=\th_s-(\s'+\a_{i_{k}})\in \D$, and 
$\s'+\a_{i_{k}}\leq \s$, against the 
maximality of $\s'$. It follows that $(\t, \a_{i_{k}})<0$, 
a contradiction. Hence $\t=\th_s'+ \u$.  By 5.1, $\th_s'$ is not dominant, 
hence $(\th'_s,\a_i)<0$ for some $i\in [n]$.  Since 
$\th_s'=\th_s-\s'$, $\th_s$ is dominant, and $(-\s', \a_{i_j})\geq 0$
for $j\in [n]\meno [h]$, we obtain that $(\th_s', \a_{i_k})<0$
for some $k\in [h]$. But then $(\u, \a_{i_{k}})\leq 0$, 
hence $(\t, \a_{i_{k}})<0$, against the assumption.
\endemo
\bigskip
 
\ni
{\bf Notation.} 
We set
$$\D_\th= \{\a\in \D\mid (\a, \th)=0\}, 
\qquad
\Pi_\th'=\{\a\in\Pi\mid (\a,\th)\not=0\}.$$
\smallskip

\proclaim{Lemma 5.3}
\roster
\item
$\D_\th$ is the standard parabolic subsystem of $\D$ generated by
$\Pi\meno\Pi'_\th$.
\item
The longest element $w_0$ of $W$ preserves the irreducible
components of $\D_\th$ and acts on them as their longest element.
\item
Let $\Gamma$ be any irreducible component of $\D_\th$, $\th(\Gamma)$ 
its highest root and $\a\in \Pi'_\th$. Then $(\a, \th(\Gamma))< 0$.
\endroster
\endproclaim

\demo{Proof}
\parni
(1) This is well known, since $\th$ is dominant.
\parni
(2) Let $w_0^\th$ be the longest element of $W(\D_\th)$. Then
$w_0^\th$ is the commutative product of the longest elements relative
to each irreducible component of $\D_\th$, hence it preserves $\D_\th$
componentwise and acts on each component as its longest element.
Consider $s_\th w_0^\th$. For all $\a\in \Dp$ we have $s_\th(\a)=
\a-(\a,\thvee)\th$, hence $s_\th(\a)<0$ if and only if $\a\not\perp\th$; 
it follows that  $N(s_\th)=\Dp\meno \Dp_\th$. Moreover, $s_\th$ fixes
$\D_\th$ 
pointwise, thus we have $N(s_\th w_0^\th)=
N(s_\th)+s_\th N(w_0^\th)=(\Dp\meno \Dp_\th)+s_\th \Dp_\th=
(\Dp\meno \Dp_\th)+\Dp_\th=\Dp$. It follows that $s_\th w_0^\th=w_0$, 
which clearly implies the claim.
\parni
(3) This follows from the fact that if $\a\in\Pi'_\th$, then $\a$ is
connected to each component of the Dynkin graph of $\Pi\meno
\Pi'_\th$. This implies that $(\a, \b)\leq 0$ for all $\b\in \Gamma$
and $(\a, \b)<0$ for at least one $\b\in \Gamma\cap \Pi$. Since
$\supp(\th(\Gamma))= \Gamma\cap\Pi$, we obtain $(\a,
\th(\Gamma))<0$.
\endemo
\bigskip
\proclaim{Definition 5.4}
Let $$\D=\D_1\sps\D_2\sps\ldots\sps \D_k$$ be a chain of standard 
parabolic irreducible subsystems of $\D$ and set $\th_i=\th(\D_i)$ 
for $i=1, \dots, k$. We shall say that $\D=\D_1\sps\D_2\sps\ldots\sps 
\D_k$ is a {\it good chain} if:
\roster
\item $\D_i$ is an irreducible component of $\{\a\in\D_{i-1}\mid
(\a,\theta_{i-1})=0\}$.
\item If $\a\in\Pi$ and $(\a, \th_{i})>0$, then $\a\perp \D_j$ for all
$j\geq i+2$.
\item  $\D_i$ contains at least a long root.
\endroster
\endproclaim
\smallskip

\ni
{\bf Notation.} 
For a good chain  $\{\D=\D_1\sps\D_2\sps\ldots\sps \D_k\}$  we set 
$$\th_i=\th(\D_i); \quad \Pi_i=\D_i\cap \Pi; \quad \Pi'_i=\{\a\in 
\Pi_i\mid (\a, \th_i)\not=0\},$$ for $i=1, \dots, k$.
\smallskip

\proclaim{\bf Lemma 5.5}
Let $\{\D=\D_1\sps\D_2\sps\ldots\sps \D_k\}$ be a good chain.
\roster
\item
For $i=1, \dots, k$, $\th_i$ is a long root, in particular
$(\th_i^\vee, \b)\in \{-1, 0,1,2\}$, for all $\b\in \Dp$
and $(\th_i^\vee, \b)=2$ if and only if $\b=\th_i$.
\item
We have $(\th_i^\vee,\b)=0$ for all $\b\in\D_j$ with $1\leq i<j\leq 
k$.
\item
Let $\a\in \Pi$. Then $(\th_i^\vee,\a)=-1$ if and only if $\a\in 
\Pi_{i-1}'$, and $(\th_i^\vee,\a)>0$ if and only if $\a\in \Pi_i'$.
\endroster
\endproclaim

\demo{Proof}
The first statement follows from property (3) in Definition 5.4. The
second statement follows directly from the definitions, too. We prove
(3). Assume $\a\in \Pi$. By Lemma 5.3 (3), if $\a\in \Pi_{i-1}'$, we
have $(\th_i^\vee,\a)=-1$. Moreover, it is clear that
$(\th_i^\vee,\a)\in \{1, 2\}$ if $\a\in \Pi_i'$ and that
$(\th_i^\vee,\a)=0$ if $\a\in \Pi_i\meno \Pi_i'$, so it remains to
prove that $(\a, \th_i)=0$ if $\a\not\in\Pi_{i-1}'\cup \Pi_i$. Assume
$\a\in \Pi_{i-j}\meno (\Pi_{i-1}'\cup \Pi_i)$, with $j\geq 1$. If
$\a\in \Pi_{i-j}'$, then $j>1$ and we are done by property (2) of
Definition 5.4. If $\a\not\in \Pi_{i-j}'$, then $\a$ lies in a
connected component of $\D_{\th_{i-j}}$. If such
component is not $\D_{i-j+1}$, then $\a\perp \Pi_{i-j+1}$, hence
$(\a, \th_i)=0$, since $\th_i\in \D_{i-j+1}$. If $\a\in\Pi_{i-j+1}$,
we conclude by induction.
\endemo
\bigskip
Set $X^+=X\meno\{0\}$.\bigskip
\proclaim{Proposition 5.6}
Let $\{\D=\D_1\sps\cdots\sps\D_k\}$ be a good chain of root subsystems
of $\D$ and set $\t=\sum\limits_{i=1}^k \th^\vee_i$. Then $\t\in X^+$.
Moreover, if $\a\in \Pi$, we have $(\t,\a)\not=0$ if and only if
$\a\in \Pi_k'$.
\par                
Conversely, for each $\t\in X^+$ there exists a good chain
$\{\D=\D_1\sps\cdots\sps\D_k\}$ such that
$\t=\sum\limits_{i=1}^k\th^\vee_i$.
\endproclaim

\demo{Proof}
Clearly $\sum\limits_{i=1}^k\th^\vee_i \in\Q$, so we only have to 
prove that $\sum\limits_{i=1}^k\th^\vee_i \in\overline C_2$. We have 
$\(\sum\limits_{i=1}^k\th^\vee_i,\th\)=(\th^\vee,\th)+\sum_{i=2}^k 
(\th^\vee_i,\th)=2$ since, for $i\geq 2$, $\th_i\in\D_\th$. Thus it 
remains to prove that $\sum\limits_{i=1}^k\th^\vee_i$ is dominant. 
Let $\a\in\Pi$. By Lemma 5.5, $(\a, \th_i)\not=0$ if and only if either 
$\a\in \Pi_i'$ or $\a\in \Pi_{i-1}'$. Set $i^*=\max\{i\mid\a\in 
\Pi_i\}.$ Then $\a\not\in \Pi_i'$ for all $i<i^*$ and $i>i^*$. If 
also $\a\not \in \Pi_{i^*}'$, we obtain 
$\(\a,\sum\limits_{i=1}^k\th^\vee_i\)=0$. Assume $\a\in \Pi_{i^*}'$.
If $i^*<k$, then clearly $\a\not=\th_{i^*}$ and by Lemma 5.5 we obtain 
$\(\a,\sum\limits_{i=1}^k\th^\vee_i\) 
=(\a,\th^\vee_{i^*})+(\a,\th^\vee_{i^*+1})= 1-1=0$. If $i^*=k$, then 
$\(\a,\sum\limits_{i=1}^k\th^\vee_i\) =(\a,\th^\vee_{i^*})\in 
\{1,2\}$.
\par
Now we prove the converse statement. Assume $\t\in X^+$. As
remarked in 4.1 we have $(\t,\th)=2$ and, by Lemma 5.2, $\t\geq
\th^\vee$. If $\t=\th^\vee$ we are done, so we assume
$\t-\th^\vee\not=0$. We have $(\t-\th^\vee, \th)=0$ and
$\t-\th^\vee\in (\Q)^+$. Since $\th$ is dominant, this implies that
all the simple coroots which support $\t-\th^\vee$ must be
perpendicular to $\th$, so that, indeed, $\t-\th^\vee\in
\Q(\D_{\th})$. It is clear that $\t-\th^\vee$ is still dominant
relatively to $\D_{\th}$; we claim that, moreover,
$\supp(\t-\th^\vee)$ is included in a single irreducible component of
$\D_\th$. By Lemma 5.3 (3), for any irreducible component $\Gamma$ of
$\D_\th$ and $\a\in \Pi'_\th$ we have $(\a, \th(\Gamma)^\vee)=-1$. Let
$\Gamma_1,\,\ldots,\Gamma_j$ be irreducible components of $\D_\th$.
Then $\t-\thvee$ is a sum of dominant elements in
$\Q(\Gamma_1),\ldots, \Q(\Gamma_j)$.
If, by contradiction, $\supp(\t-\th^\vee)$ has non-trivial
intersection with two distinct irreducible components, say $\Gamma_1$
and $\Gamma_2$, of $\D_\th$, then, by Lemma 5.2, $\t-\th^\vee \geq
\th(\Gamma_1)^\vee+\th(\Gamma_2)^\vee$. It follows that, for
$\a\in\Pi'_\th$, we have  $(\a, \t-\th^\vee)\leq -2$, and therefore
that $(\a, \t)\leq -1$, which is impossible. So there exists an
irreducible component $\D_2$ of $\D$ such that $\t-\th^\vee$ lies in
$\Q(\D_2)$. By Lemma 5.2, $\t-\th^\vee\geq \th_2^\vee$. By induction
we can find a chain of irreducible subsystems $\D=\D_1\sps\D_2\sps
\cdots \sps \D_k$, such that $\D_i$ is an irreducible component of
$(\D_i)_{\th_i}$ and $\tau=\sum\limits_{i=1}^k\th_i^\vee$. 
\par
It remains to prove that $\D=\D_1\sps\D_2\sps \cdots \sps \D_k$
satisfies conditions (2) and (3) of Definition 5.4. 
\par
First we prove (3). Assume by contradiction that $\D_i$ consists
entirely of short roots, for some $i$ in $[k]$. By a direct check we
obtain that then $\D_i$ has rank $1$ and $\D_{i-1}$ is of type $B_3$
or $G_2$; in particular $i=k$. Now $\th_{k-1}^\vee+\th_k^\vee$ must be
dominant for $\D_{k-1}$, but this is not the case. In fact, $\th_k$ is
the short root of $\Pi_{k-1}$, while $\th_{k-1}$ is long. If $\a$ is the (unique) root in
$\Pi'_{k-1}$, then $\a$ is long, thus $(\th_{k-1}^\vee,\a)=1$ and $(\th_k^\vee,\a)<-1$, 
hence $(\th_{k-1}^\vee+\th_k^\vee,\a)<0$.
\par
Finally we prove (2). Assume $\a\in \Pi_i'$. Clearly $(\a,\th_j)=0$
for $1\leq j<i$. Moreover, as remarked above, we have $(\a,
\th_{i+1})<0$ and, since $\a\not\in\D_j$ for $j>i$, $(\a, \th_j)\leq
0$ for $j\geq i+2$. Since $\th_i$ and $\th_{i+1}$ are long we obtain
$(\t,\a)=(\th_i^\vee,\a)+(\th_{i+1}^\vee,\a)+\sum_{j\geq
i+2}(\th_j^\vee,\a)= 1-1+\sum_{j\geq i+2}(\th_j^\vee,\a)$. Since $\t$
is dominant it follows that $(\th_j^\vee,\a)= 0$ for all $j\geq i+2$.
\endemo
\bigskip

\ni
{\bf Remark 5.7.}
The proof of Theorem 5.6 makes clear the following facts.
\roster 
\item
Assume $\t\in X^+$ and $\t=\th_1^\vee+\cdots+\th_k^\vee$, as in 5.6. 
Then ${(\th, \th)\over 2}\t=\th_1+\cdots+\th_k$, and ${(\th, \th)\over 2}
(\t,\t)=2k$.
\item 
Assume $\t,\tilde\t\in X^+$, $\t=\th_1^\vee+\cdots+\th_k^\vee$, and
$\tilde \t=\tilde\th_1^\vee+\cdots+\tilde\th_{\tilde k}^\vee$ 
as in Theorem 5.6.
Let $h$ be the maximal index such that $\th_h=\tilde\th_h$. Then
$\th_i=\tilde\th_i$ for $1\leq i\leq h$; moreover,
$\th_i\perp\tilde\th_j$ for all $i, j>h$. It follows that ${(\th,
\th)\over 2} (\t,\tilde \t)=2 h$.
In particular the orthogonal decomposition of $\t\in X^+$ provided by
Theorem 5.6 is unique.
\item 
Let $\a\in\Dp_\ell$. Then there exists a unique maximal good chain 
$\D=\{\D_1\sps\cdots\sps \D_k\}$ such that $\a\in \D_i$ for
$1\leq i\leq k$. Moreover, if $\{\D=\D'_1\sps\cdots\sps \D'_h\}$
is any good chain such that $\a\in \D'_h$, then $h\leq k$ and $\D'_i=\D_i$ for 
$1\leq i\leq h$. This is easily seen by induction; indeed, if 
$\{\D=\D_1\sps\cdots\sps \D_i\}$ is a good chain such that $\a\in \D_i$, 
then $\a$ belongs to at most one connected component
of $\D_i\cap \th_i^\perp$, hence there is at most one good chain 
$\{\D=\D_1\sps\cdots\sps \D_{i+1}\}$ which extends 
$\{\D=\D_1\sps\cdots\sps \D_i\}$ and such that $\a\in \D_{i+1}$.
\endroster
\bigskip

\ni 
{\bf Examples.}
\smallskip

\ni
$\D\cong A_{n}$. Set $\Delta_i=\langle \a_i,\ldots,\a_{2n-i+1}\rangle\
1\leq i\leq \lfloor \frac{n+1}{2}\rfloor$. Then the possible good
chains are $\Cal S_j=\{\D_1,\ldots, \D_j\},\ 1\leq j\leq \lfloor
\frac{n+1}{2}\rfloor$. The corresponding elements in $X$ are
$\t_j=\omega^\vee_j+\omega^\vee_{n-j+1}$ (here $\omega^\vee_i$ denotes
the $i$th fundamental coweight).
\smallskip

\ni
$\D\cong E_6$. Due to condition (2), there are only two good chains:
$$\{\D\},\qquad\{\D,\langle \a_1, \a_3, \a_4, \a_5, \a_6\rangle\}.$$
Indeed $X^+=\{\thvee,\,\omega_1^\vee+\omega_6^\vee\}.$
\smallskip

\ni
$\D\cong B_{2n+1}$. Due to conditions (2) and (3), there are the
following good chains: $$\{\D\},\quad\{\D,\langle\a_3, \a_4,\dots,
\a_{2n+1}\rangle\},\ldots$$ $$\{\D,\langle\a_3, \a_4,\dots,
\a_{2n+1}\rangle,\ldots ,\langle\a_{2n-1}, \a_{2n},\dots,
\a_{2n+1}\rangle\},\ \{\D,\langle \a_1\rangle\}.$$
\vskip 1cm

\heading \S6 
The posets $\II(\t)$
\endheading
\medskip

In this section we begin to study  the poset of  abelian ideals encoded by
elements in $\widehat W$ of the form $t_\t v$, for a fixed $\t\in X$
and $v$ varying in $W$. We prove that this poset has minimum, obtained
for a special element $t_\t v_\t$ which will be analysed in detail. Then we introduce the posets
$\II(\a,\t)$, which will be extensively studied in the next section.
\bigskip

\ni
{\bf Notation.} Let $\tau\in X^+$. According to Proposition 5.6, and
Remark 5.7 (2) there exists a unique good chain of irreducible
subsystems $\{\D=\D_1\sps \cdots\sps \D_k\}$ such that $\tau =
\th^{\vee}_1+\cdots+\th^{\vee}_k$, where $\th_i$ is the highest root of
$\D_i$. We shall call $\tau = \th^{\vee}_1+\cdots+\th^{\vee}_k$ {\it the
standard orthogonal decomposition of} $\t$, and $\{\D=\D_1\sps
\cdots\sps \D_k\}$ {\it the good chain of $\t$}. Moreover we set 
$$\D(\tau)=\D_k,\quad \Pi(\tau)=\Pi\cap\D(\tau), \quad
\th(\tau)=\th_k, \quad \Pi'(\tau)=\{\a\in\Pi(\t)\mid
(\a,\th(\t))\not=0\}.$$
We set $\D(0)=\Pi(0)=\emptyset$. Then, by 5.6, given $\a\in \Pi$   we
have $(\a,\t)\not=0$ if and only if $\a\in \Pi'(\t)$.
\par
\bigskip

We set, for $\t\in X$,
$$\wab(\t)=\{w\in\wab\mid  w=t_\t v,\ v\in W\},$$ 
$$\II(\t)=\{\i\in\II\mid w_\i\in\wab(\t)\}.$$

\ni
{\bf 6.1.} 
Let $w\in \hatw$. 
By Proposition 3.1 (2), $w\in\wab$ if and only if $w(C_1)\sbs C_2$.
Let $x\in C_1$, and $w=t_\t v$, with $\t\in \Q$ and $v\in W$. Then
$w(x)\in C_2$ if and only if $(t_\t v(x), \b)>0$ for all $\b\in \Pi$,
and $(t_\t v(x), \th)<2$. Equivalently, $w(x)\in C_2$ if and only if
$0<(t_\t v(x), \b)<2$ for all $\b\in \Dp$. For $\b\in \D$ have $(t_\t
v(x), \b)=(\t+v(x), \b)= (\t,\b)+(x, v^{-1}(\b))$; moreover $0<|(x,
v^{-1}(\b))| <1$ and $(x, v^{-1}(\b))>0$ if and only if
$v^{-1}(\b)>0$. Set
$$\D^i_\t=\{\a\in \D\mid (\a, \t)=i\}.$$ Then we obtain:  
\bigskip
\item{(1)}
$t_\t v\in \wab$ if and only if the following conditions hold: 
$\t\in X$; for all $\b\in \D_\t^0\cap \Dp$, $v^{-1}(\b)>0$; 
for all $\b\in \D_\t^2$, $v^{-1}(\b)<0$.
\bigskip
By 4.1, for $\t\in X$ we have $(\t, \th)=2$ unless $\t=0$. In this case, 
by 3.1 (2), we have that  $t_\t v\in \wab$ if and only if $t_\t v=1$, the 
identity of $\hatw$. 
So we have, equivalently, 
\bigskip
\item{(2)} 
$t_\t v\in \wab$ if and only if the following conditions hold: 
$\t\in X$; for all $\b\in \Pi$, if $(\t, \b)=0$,  then $v^{-1}(\b)>0$; 
if $t_\t v\ne 1$, then $v^{-1}(\th)<0$.
\bigskip
Using Theorem 5.6, we also obtain:
\bigskip 
\item{(3)}
$t_\t v\in \wab$ if and only if $\t\in X,\ L_v\subseteq \Pi'(\t)$, and, 
if $t_\t v\ne 1$, $\th\in N(v)$.
\bigskip
Moreover, by Proposition 3.3,  $w\in \wab$ if and only if
$N(w)\sbsq\Dm+\d$. We have
$w^{-1}(-\b+\d)=v^{-1}t_{-\t}(-\b+\d)=v^{-1}(-\b+(1-(\t,\b))\d)=
v^{-1}(-\b)+(1-(\t,\b))\d$. Hence $w^{-1}(-\b+\d)<0$ if and only if
$(\t,\b)>1$ or $(\t,\b)=1$ and $v^{-1}(\b)>0$. 
Thus if $w\in \wab$, since $(\t, \b)\leq 2$ for all $\t\in X$ and
$\b\in \Dp$, we obtain (cfr.\cite{CMP, Lemma C}):
\bigskip
\item{(4)}
$N(t_\t v)=\{-\b+\d\mid \b\in \D^2_\t\cup(\D_\t^1\meno N(v))\}$,  
\bigskip
hence:
\bigskip
\item{(5)}
$\gog_\b\sbsq \goi(t_\t v) \text{ \ if and only if \ \ }\b\in
\D^2_\t\cup(\D_\t^1\meno N(v))$. 
\bigskip
By (1), if $w\in \wab$, then $\D^2_\t\sbsq N(v)$. Therefore by (5),  
if $t_\t v$ and $t_\t v'$ belong
to $\wab$ ($\t\in X$, $v, v'\in W$), then the inclusion of the corresponding ideals corresponds
to the reverse weak order of $v,\, v'$: 
\bigskip
\item{(6)} 
$\goi(t_\t v)\sbsq \goi(t_\t v') \text{ if and only if } v'\weak v$. 
\bigskip
Consider
$$\D_\t^{even}=\{\a\in \D\mid (\a, \t)\text{ is even}\}.$$ 
It is clear that for $\t\in X$
\bigskip
\item{(7)} 
$\D^2_\t=(\D_\t^{even}\meno\D_\t^0)\cap \Dp$.
\bigskip 
We remark that $\D_\t^0$ is the parabolic standard subsystem 
of $\Da$ generated by  $\Pi\meno \Pi'(\t)$, and,
by \cite{CMP, Lemma 2.6 (2)}, $\D_\t^{even}$ is isomorphic to the standard
parabolic subsystem of $\Da$ generated by $\Pia\meno \Pi'(\t)$. 
\bigskip

\proclaim{Proposition 6.2}
Let $w\in \wab$, $\t\in X^+$, and let $\t=\th_1^\vee+\cdots+\th_k^\vee$
be the standard orthogonal decomposition of $\t$. 
Then $\goi(w)\in\II(\t)$ if and only if $-\th_i+\d\in N(w)$ for all
$i\in \{1, \dots, k\}$ and $-\b+\d\not\in N(w)$ for all $\b\in \D_k
\cap\th_k^\perp$.
\endproclaim

\demo{Proof}
Assume $w=t_\t v$ with $v\in W$. 
Since $(\t, \th_i)=2$ for all $i\in \{1, \dots, k\}$, by 6.1 (4) we have
that $w^{-1}(-\th_i+\d)<0$, for all $i\in \{1, \dots, k\}$.
Moreover, if $\b\in \D_k\cap\th_k^\perp$, we have $(\t,\b)=0$, 
hence by 6.1 (4) $-\b+\d\not\in N(w)$. 
Conversely, let $w=t_{\tilde\t} v$ with $\tilde\t\not=\t$ and
$\tilde\t=\tilde\th_1^\vee+\cdots +\tilde\th_{\tilde k}^\vee$ be the
standard orthogonal decomposition of $\tilde\t$. 
Let $h\in \{1,\dots, \min\{k,\tilde k\}\}$ be the maximal index such
that $\th_h=\tilde\th_h$.  Then, by Remark 5.7 (2) at least one of
$k$, $\tilde k$ is strictly greater than $h$. Moreover, if $h<k$, then
$\th_{h+1} \perp \tilde \t$ hence, by 6.1 (4), $-\th_{h+1}+\d\not\in N(w)$. If
$h= k$, then we have that $\tilde\th_{h+1}\in\D_k
\cap\th_k^\perp$ and  $-\tilde\th_{h+1}+\d
\in N(w)$.
\endemo
\bigskip

As recalled in 2.9,  for each $I\sbsq \Pi$ there is a unique $v\in W$
with $L_v\sbsq I$ of maximal length, and maximal, indeed, with respect
to the weak order: just take $v=w^I_0 w_0$, where $w_0^I$ is the
longest element in $\la s_\a\mid \a\not\in I\ra$. Then
$N(v)=\Dp\meno\D^I$, where $\D^I$ is the subsystem of $\D$
generated by $\Pi\meno I$.
\bigskip

\ni
{\bf Notation.} 
We fix $\t\in X$ and set  $$v_\t=w^I_0 w_0, \qquad \text{\rm with} \qquad
I=\Pi'(\t).$$ 
\bigskip

By 6.1 (3) $t_\t v_\t\in \wab$. Moreover we have the following result.
\bigskip

\proclaim{Proposition 6.3} 
\roster 
\item
$\goi(t_\t v_\t)=\min\II(\t)$;
\item 
$\dim\goi(t_\t v_\t)={1\over 2}(|\D_\t^{even}|-|\D_\t^0|)$.
\endroster
\endproclaim

\demo{Proof} (1) From 6.1 (6) we directly obtain that $\goi(t_\t
v_\t)$ is minimal in $\II(\t)$.
\parni (2) We have $N(v_\t)= \Dp\meno \D^I=\Dp\meno
\D_\t^0=\D_\t^2\cup \D_\t^1$, hence, by 6.1 (5), $\goi(t_\t
v_\t)={\mathop\oplus\limits_{\a\in \D_\t^2}}\gog_\a$. In particular
$\dim\goi(t_\t v_\t)=|\D_\t^2|$, and by 6.1 (7) we obtain formula (2).
\endemo
\bigskip

\proclaim{Lemma 6.4}
The element $v_\t$ is an involution. 
\endproclaim

\demo{Proof}
We have, for any $I\sbsq \Pi$,  that $(w_0^I w_0)^{-1}=w_0^J
w_0$ with $J=-w_0(I)$. So, in our case,  we have to prove that
$-w_0(I)=I$, for $I=\Pi'(\t)$. By Lemma 5.3 (2) $w_0$ preserves all
irreducible components of $\D_\th$ and acts on them as their longest
element. Since the longest element maps the highest root to its opposite, 
by recursion we obtain that $-w_0$ fixes $\t$, for each
$\t\in X$, and therefore stabilizes $\Pi'(\t)$.
\endemo
\bigskip

\proclaim{Lemma 6.5}
Let $I\sbsq \Pi$, $v=w^I_0w_0$, and $y\in W$. If $v$ is an involution,
then $L_{vy}\sbsq I$ if and only if $L_y\sbsq I$. In this case
$-vN(y)\sbsq N(v)$ and  $N(vy)=N(v)\meno -vN(y)$. In particular
$l(vy)=l(v)-l(y)$.
\endproclaim

\demo{Proof} 
For any $u\in W$ we have $L_{u}\sbsq I$ if and only if $\Npm(u)\sbsq
\D\meno \D^I$, where $\D^I$ is the subsystem of $\D$ generated by
$\Pi\meno I$. Recall that $\Npm(vy)=\Npm(v)+v\Npm(y)$.  If
$L_y\not\sbsq I$, let $\a\in L_y\meno I$. Then $\a\in N(y)$,
whence $v(\a)\in vN(y)$, but $v(\a)=v^{-1}(\a)>0$, hence $v(\a)\not\in
\Npm(v)$; therefore $v(\a)\in \Npm(vy)$. Since $\Npm(v) = \D\meno
\D^I$, we obtain that $\Npm(vy)\not\sbsq \D\meno \D^I$, hence
$L_{vy}\not\sbsq I$. Conversely, if $\Npm(y)\sbsq \D\meno \D^I$,
then $vN(y)=v^{-1}N(y) <0$. By 2.7, it follows that $N(vy)\sbsq N(v)$,
hence $L_{vy}\sbsq L_v$. Moreover, $-v(N(y))\sbsq N(v)$ and
$l(vy)=l(v)-l(y)$.
\endemo
\bigskip

\proclaim{Proposition 6.6}
Assume $\t\in X$, $v, v'\in W$. Then $t_\t v\in \wab$ if and only if 
$v=v_\t y$ with $L_y\sbsq \Pi'(\t)$. 
\par
Assume $t_\t v,\ t_\t v'\in \wab$, $v=v_\t y$, and  $v'=v_\t y'$. 
Then $\i(t_\t v)\sbsq \i(t_\t v')$ if and only if $N(y)\sbsq N(y')$.
\endproclaim

\demo{Proof}
The first assertion follows directly from 6.1 (3), 6.4, and 6.5. 
\par
Under the assumptions of the second assertion, by 6.5 we have 
$-v_\t N(y)\sbsq N(v_\t)$, $-v_\t N(y')\sbsq N(v_\t)$, 
$N(v_\t y)=N(v_\t)\meno -v_\t N(y)$, and 
$N(v_\t y')=N(v_\t)\meno -v_\t N(y')$. 
Therefore, $v'\weak v$ if and only if $y\weak y'$. This proves the claim.
\endemo
\bigskip

It is clear that the arguments of section 4 can be specialized
to $\II(\t)$. We obtain the following result.
\bigskip

\proclaim{Proposition 6.7}
Let $\goi\in \II(\t)$. Then $\goi$ is maximal in $\II(\t)$ if and only if 
$w_\goi(\Pi)\cap (\Dm+\d)=\emptyset$. In this case, if $\t\ne 0$, we have 
$w_\goi^{-1}(-\th+2\d)\in \Pi$, hence if $w_\goi=t_\t v$, then $v^{-1}(-\th)
\in \Pi_\ell$. 
\endproclaim

\demo{Proof}
We notice that if $\goj\in \II(\t)$, then $w_\goj=w_\goi x$ with $x\in
W$. So,  we can just mimic the proof of 4.4 considering $\Pi$ in
place of $\Pia$.
\endemo
\bigskip

The above result leads us to analyse $v^{-1}(\th)$ for all $v\in W$ such 
that $t_\t v\in \II(\t)$. We first study what $v_\t(\th)$ is. 
\bigskip

For any root $\b$ and any simple root $\a$ let $n_\a(\b)$ denote the 
coefficient of $\a$ in $\b$ and set $m_\a=n_\a(\th)$: 
$$\b=\sum_{\a\in\Pi}n_\a(\b)\a,\qquad\th=\sum_{\a\in\Pi}m_\a\a.$$
\smallskip

\proclaim{Proposition 6.8} Fix $\t\in X^+$.
\roster 
\item We have
$v_\t(\Pi\meno \Pi'(\t))=\Pi\meno \Pi'(\t)$. 
\item  
Set $B_\t=v_\t(\Pi'(\t))$. For each $\b\in B_\t$, there exists a 
unique $\a\in \Pi'(\t)$ such that $n_\a(\b)\not=0$ and, for such an 
$\a$, we have $n_\a(\b)=-1$. In particular $B_\t\sbsq \Dm$.
\item  
Let $\b\in B_\t$ and $\a\in \Pi'(\t)$ such that $n_\a(\b)=-1$. Then 
$-\b$ is the maximal element in $\{\g\in \Dp\mid 
\supp(\g)\cap\Pi'(\t)= \{\a\},\ n_\a(\g)=1\}$, with respect to the 
standard partial order of $Q$.
\item 
We have $v_\t(\th)=-\th(\t)$. Moreover, $n_\a(\th(\t))=m_\a$ for 
all $\a\in \Pi'(\t)$ and $\th(\tau)$ is the minimal root with this 
property, with respect to the standard partial order of $Q$. 
\endroster
\endproclaim

\demo{Proof}
\parni (1) 
The claim follows since $v_\t$ maps $\Pi\meno \Pi'(\t)$ to 
$-w_0(\Pi\meno \Pi'(\t))$, which is equal to $\Pi\meno 
\Pi'(\t)$, as seen in the proof of Lemma 6.4.
\par\ni 
(2) 
Recall that $v_\t=v_\t^{-1}$; thus $B_\t=v_\t(\Pi'(\t))$ is a set of
negative roots. By (1) $v_\t$ preserves $\Pi\meno \Pi'(\t)$, and since
$v_\t(\Pi)$ is a basis of $\D$, we obtain that $B_\t\cup(\Pi\meno
\Pi'(\t))$ is a basis of $\D$. In particular, each $\a\in \Pi'(\t)$,
is a nonnegative or nonpositive linear combination of
$B_\t\cup(\Pi\meno \Pi'(\t))$. This implies that for each $\a\in
\Pi'(\t)$, some $\b\in B_\t$ is supported by $\a$ and, moreover, $\a$
has coefficient $-1$ in such a $\b$. Moreover, each $\b\in B_\t$ is
supported by a unique $\a\in \Pi'(\t)$.
\parni 
(3)  
Let $\a\in \Pi'(\t)$ and $\b_\a$ be the unique root in $B_\t$
supported by $\a$. If $\g\in \Dp$, $\supp(\g)\cap \Pi'(\t)=\{\a\}$,
and $n_\a(\g)=1$ then, arguing as in (2) we find that $\g$ equals
$-\b_\a$ minus a sum of roots in $\Pi\meno \Pi'(\t)$. This proves the
claim.
\par\ni (4) 
We have $v_\t(\th)=\sum_{\a\in\Pi}m_\a v_\t (\a)$. For each  $\a\in
\Pi'(\t)$, set $\g_\a=v_\t(\a)$, and let $\a'$ be the unique root of
$\Pi'(\t)$ which supports $\g_\a$.
Then, by (2), $n_{\a'}(v_\t(\th))=-m_\a$, hence $m_\a\leq m_{\a'}$.
Moreover $\{\a'\mid \a\in \Pi'(\t)\}= \Pi'(\t)$, hence $\{m_{\a'}\mid
\a\in \Pi'(\t)\}=\{m_\a\mid \a\in\Pi'(\t)\}$. It follows that
$\sum\limits_{\a\in\Pi'(\t)} m_\a= \sum \limits_{\a\in\Pi'(\t)}
m_{\a'}$, hence $m_\a= m_{\a'}$.
Thus each $\a\in \Pi'(\t)$ occurs in $-v_\t(\th)$ with coefficient
$m_\a$. We show that $-v(\th_\t)$ is the minimal root with this property.
As seen in (1), $v_\t$ permutes $\Pi\meno \Pi'(\t)$; since $\th+\eta$
is not a root for each $\eta\in Q^+$, we obtain that $v_\t(\th)+\eta$
is not a root for all $\eta\in Q^+(\Pi\meno \Pi'(\t))$. Equivalently
$-v_\t(\th)-\eta$ is not a root for all $\eta\in Q^+(\Pi\meno
\Pi'(\t))$, hence $-v_\t(\th)$ is minimal among the roots in which
$\a$ occurs with coefficient $m_\a$ for all $\a\in \Pi'(\t)$. Now we
know from the orthogonal decomposition of
$\t=\thvee_1+\cdots+\thvee_k$ that $(\t,\th_i)=2$ for $i=1,\ldots,k$.
Since $n_\a(\th_i)\leq n_\a(\th_1)=m_\a$ for all $\a\in \Pi'(\t)$ and
since $\Pi'(\t)$ is the set of all simple roots non orthogonal to
$\t$, we obtain $n_\a(\th_i)=m_\a$ for $i=1,\ldots,k$. In particular
we obtain that $-v_\t(\th)\leq \th(\t)$. But $\Pi'(\t)$ is the set of
all simple roots which belong to $\D(\t)$ and are not orthogonal to
$\th(\t)$, hence we have $(-v_\t(\th),
\thvee(\t))=(\th(\t),\thvee(\t))=2$. Since $\th(\t)$ is long this
implies that $-v_\t(\th)= \th(\t)$.
\endemo 
\bigskip

\proclaim{Lemma 6.9} 
Let $\a\in\Dp_\ell,\,v\in W$ and assume that $\a\notin\D(\tau)$ and 
$v(\a)=\th(\tau)$. Then $L_v\not\sbsq \Pi'(\t)$.
\endproclaim 

\demo{Proof} 
By 2.3 (5), $L_v\not\sbsq \Pi'(\t)$, if and only if there exists
some $\g\in N(v)$ such that $\supp(\g)\cap\Pi'(\t)= \emptyset$. We
shall find such a $\g$. Since $\a\not\in\D(\tau)$, we have that $v\notin
W(\D(\tau))$, hence $N(v)\not\sbsq \D(\tau)$.  By  2.8, $N(v)\cap\D(\tau)=
N(v')$ for
some $v'\in W(\D(\tau))$, hence there exists $v''\in W$ such that
$v=v'v''$, and $v'N(v'')=N(v)\meno \D(\tau)$. Set $N(v'')=\{\b_1,
\dots, \b_h\}$, and $\g_i=v'(\b_i)$, for $i=1,\dots, h$, so that  
$N(v)\meno \D(\tau)=\{\g_1, \dots, \g_h\}$. By 2.4 (1), 
we may assume that $v=s_{\g_h}\cdots s_{\g_1} v'$. Conditions $\a\notin
\D(\t)$ and $v'\in W(\D(\t))$ imply $v'(\a)\ne\th(\t)$, hence
$s_{\g_h}\cdots s_{\g_1}(\th(\t))\ne \th(\t)$.
This implies that there exists $i\in\{1, \dots, h\}$ such that
$\g_i\not\perp \th(\t)$. Since $\th(\t)\notin N(v)$, we have that
$\th(\t)\ne \g_i$, hence $(\g_i, \thvee(\t))=\pm 1$.
Assume first $(\g_i, \thvee(\t))=-1$, so that $\g_i+\th(\t)\in \Dp$.
By 6.8 (4), for any $\b\in \Pi'(\t)$ we have that
$n_\b(\th(\t))=m_\b$, hence that $n_\b(\g_i)=0$, and we are done. So
assume $(\g_i, \thvee(\t))=1$, so that $\g_i-\th(\t)\in\D$. Since
$\g_i\notin\D(\t)$, we have that $\g_i\not<\th(\t)$, hence
$\g_i=\th(\t)+\g_i'$, for some $\g_i'\in \Dp$. Since $\g_i\in N(v)$,
and $\th(\t) \notin N(v)$, by 2.3 (3) we obtain that $\g'_i\in N(v)$.
But $(\g'_i, \thvee(\t))=-1$, so we are reduced to the previous case.
\endemo
\bigskip

For $\t\in X$ and $\a\in \Dp_\ell$, we set
$$\wab(\t,\a)=\{w\in\wab(\t)\mid w(\a)=-\th+2\d\};$$
$$\II(\t,\a)=\{\goi\in\II\mid w_\goi\in\wab(\t,\a)\}.$$ 
Let $\a\in\Dp_\ell$ and $w\in\wab$ be such that $w(\a)=-\th+2\d$.
Then, if $w=t_\t v$ with $\t\in\Q$ and $v\in W$, we have $v(\a)=-\th$.
Conversely, if $v(\a)=-\th$, then $v\ne 1$, hence by 3.1 (2) $\t\ne
0$. By 4.1 it follows that $(\t, \thvee)=2$, hence $t_\t
v(\a)=t_\t(-\th)=-\th-(\t,-\th)\d= -\th+ 2\d$. Therefore
$$\wab(\t,\a)=\{t_\t v \in \wab(\t)\mid v(\a)=-\th\}.$$
It is clear that $\wab(\t, \a)$ is empty for $\t=0$. If $\t\in X^+$
and $t_\t v\in \wab$, then by 6.1 (2) we have that $v^{-1}(-\th)\in \Dp_\ell$, 
hence $$\wab(\t)=\bigcup\limits_{\a\in \Dp_\ell}\wab(\t, \a), 
\quad\text{\rm and}\quad
\II(\t)=\bigcup\limits_{\a\in \Dp_\ell}\II(\t, \a)$$
\smallskip

\proclaim{Proposition 6.10}
Assume  $\t\in X^+$  and $\a\in \Dp_\ell$.
If $\II(\t,\a)\ne\emptyset$, then $\a\in\D(\tau)$.
\endproclaim

\demo{Proof}
Assume $t_\t v\in \II(\t,\a)$. Then by  6.6 $v=v_\t y$ with
$L_y\sbsq \Pi'(\t)$; moreover, using 6.4, and 6.8 (4),
$y(\a)=v_\t^{-1}v(\a)=v_\t(-\th)=\th(\t)$. By 6.9,  this implies
that $\a\in \D(\t)$.
\endemo
\vskip 1cm
  
\medskip
\heading \S7 The posets $\II(\t,\a)$\endheading
\medskip
In this section we prove that, for $\a\in\Dp_\ell$, $\II(\t, \a)$ is
not empty if and only if $\a\in \D(\t)$, and in this case  $\II(\t,
\a)$ has a unique minimal element and a unique maximal element.
\bigskip

Let 
$$\thvee=\sum\limits_{\a\in \Pi}m^\vee_\a \a^\vee.$$ 
We have
$\thvee={2\over (\th, \th)}\th ={2\over (\th, \th)}\sum\limits_{\a\in
\Pi}m_\a {(\a, \a)\over 2}\a^\vee$, hence 
$m_\a^\vee=\frac{(\a,\a)}{(\th,\th)}m_\a$.
\par
The dual Coxeter number of $\gog$ is $$g=1+\sum\limits_{\a\in
\Pi}m^\vee_\a.$$
\smallskip

Part of next proposition is also proved in \cite{S} and \cite {P}.
\smallskip

\proclaim{Proposition 7.1}  
Let $\a$ be any long root. Then: 
\roster 
\item 
$\{v\in W\mid v(\a)=\th\}$ has a unique element $y_\a$ of 
minimal length.
\item 
$L_{y_\a}\sbsq \Pi'_\th$; conversely,
if $v'\in W$, $v'(\a)=\th$, and $L_{v}\sbsq\Pi'_\th$, then $v'=y_\a$.
\item
If $\a\in \Pi_\ell$, then $\ell(y_\a)=g-2$;
\item 
$y_\a^{-1}$ is the unique minimal element, with respect to the
weak order, in \break $\{v\in W\mid v(\th)=\a\}$; 
\item
$N(y_\a^{-1})=\{\b\in \Dp\mid (\b,\a^\vee)=-1\}$.
\endroster
\endproclaim

\demo{Proof} 
\parni
(1), (4). 
We start with proving that there exists a unique minimal
element, with respect to the weak order, in $\{v\in W\mid
v(\th)=\a\}$. This will imply (1) and (4).
Let $m=\mathop{\min}\limits_{v\in W} \{\ell(v)\mid
v(\th)=\a\}$ and let $v\in W$ be any element such that $v(\th)=\a$ and
$\ell(v)=m$.
Let $v'\in W$ be such that $v'(\th)=\a$. Then $v'=vx$ with $x\in
\stab_W(\th)$. We shall prove that $\ell(v')=\ell(v)+\ell(x)$: by
2.6, this implies in particular that $N(v)\sbsq N(v')$. Let
$v=s_{i_1}\cdots s_{i_k},\, x=s_{j_1}\cdots s_{j_h}$ be reduced
expressions. If by contradiction $s_{i_1}\cdots
s_{i_k}s_{j_1}\cdots s_{j_h}$ is not reduced, then there exist $u<_B
v,\,y<_B x$, $<_B$ being the Bruhat order, such that $uy=vx$. But
$\stab_W(\th)=W(\D_\th)=\langle s_\b \mid \b \in \Pi\meno
\Pi'_\th\rangle$, therefore, since $x\in \stab_W(\th)$, we have $y\in
\stab_W(\th)$, too. It follows that $u(\a)=\th$, against the
assumption that $\ell(v)$ is minimal. We set $y_\a=v^{-1}$; then 
it is clear that $y_\a$ is the unique element of minimal length 
in $W$ which maps $\a$ to $\th$. 
\parni
(2) 
For all $y\in W$, if  $y(\a)=\th$, then $y=x y_\a$ with 
$x\in \stab_W(\th)$ and $\ell(y)=\ell(x)+\ell(y_\a)$. In particular, 
$\ell (s_\b y_\a)>\ell(y_\a)$ for all $\b\in \Pi\meno \Pi_\th'$, hence 
$L_{y_\a}\subseteq \Pi'_\th$.
Moreover, $L_x\sbsq \Pi\meno \Pi_\th'$ and
$N(x)\sbsq N(y)$, hence, if $y\ne y_\a$, $L_{y}\not\sbsq \Pi'_\th$.
\parni
(5)
Assume $y_\a(\b)<0$. Then $\a\ne \b$ and $-y_\a(\b)\in N(y_\a)$. By item (2), \break
$\supp(-y_\a(\b))\cap \Pi'_\th\ne \emptyset$, hence 
$(\b, \a^\vee)=-(-y_\a(\b), \thvee)=-1$. Conversely, assume 
$(\b, \a^\vee)=-1$. Then $(y_\a(\b), \thvee)=-1$, hence $y_\a(\b)<0$.
\parni
(3) 
Now we assume $\a\in \Pi$ and prove that $m=g-2$. Let $\a\in \Pi_\ell$
and $v\in W$ be such that $v(\a)=\th$. We first show that $\ell(v)\geq
g-2$.  Let $v= s_{\g_1}\cdots s_{\g_k}$ be a reduced expression. Set
$v_0=1$, and $v_h=\prod_{i=1}^h s_{\g_i}$, for $1\leq h\leq k$.
Since $N(v_h)\sbsq N(v)$, $v_h^{-1}(\th)$ is a long positive root for
$0\leq h\leq k$; in particular, $v_{h-1}^{-1}(\th)\ne \g_h$, for $0<
h\leq k$.
For any  positive root $\b$ and any $\g\in \Pi$, we have
$s_\g(\b)=\b-(\b,\g^\vee)\g$. Moreover, if $\b$ is long and $\b\ne \g$, then $|(\b^\vee,\g)|\leq
1$, hence $|(\b,\g^\vee)|=\Big|(\b^\vee,\g){(\b,
\b)\over (\g, \g)}\Big|\leq \Big|{(\b, \b)\over (\g,
\g)}\Big|=\Big|{(\th, \th)\over (\g, \g)}\Big|$.  Therefore
$v_{h-1}^{-1}(\th)-v_h^{-1}(\th) \leq {(\th, \th)\over (\g_h, \g_h)}
\g_h$, for $1\leq h\leq k$, hence
$\th-\a=v_0^{-1}(\th)-v_k^{-1}(\th)\leq \sum\limits_{i=1}^k{(\th,
\th)\over (\g_i, \g_i)} \g_i$.  Since $\th=\sum\limits_{\g\in\Pi}m_\g
\g= \sum\limits_{\g\in\Pi}m^\vee_\g {(\th, \th)\over (\g,\g)}\g$, we
obtain that each $\g\ne\a$ occurs in the sequence $(\g_1, \dots,
\g_k)$ at least $m_\g^\vee$ times, and $\a$ at least
$m_{\a}^\vee-1$ times. It follows that  $\ell(v)=k\geq 
\sum\limits_{\g\in\Pi}m^\vee_\g -1=g-2$. Next, we show that there
exists a $w\in W$ such that $w(\a)=\th$ and $\ell(w)\leq g-2$. If $\b$
is a long root and $\b\not=\th$, then there exists a simple root $\g$
such that $(\b,\g^\vee)<0$ and hence $s_\g(\b)=\b+{(\th, \th)\over
(\g,\g)} \g$. Therefore, if $\a$ is a simple long root, we can find a
sequence of simple roots $(\g_1, \dots ,\g_k)$ such that
$(s_{\g_{i-1}}\cdots s_{\g_1}(\a),\g_i^\vee)<0$, $s_{\g_i}\cdots
s_{\g_1}(\a)= s_{\g_{i-1}}\cdots s_{\g_1}(\a)+ {(\th, \th)\over (\g_i,
\g_i)}\g_i$, for $i\leq k$, and $s_{\g_k}\cdots s_{\g_1}(\a)=\th$. 
Then clearly $k= \sum\limits_{\g\in\Pi}m^\vee_\g -1=g-2$ and therefore
$\ell(s_{\g_k}\cdots s_{\g_1})\leq g-2$. This concludes the proof.
\endemo
\bigskip

\proclaim{Proposition 7.2}
Let $\a\in\Dp_\ell$. Then 
for any pair $\b_1, \b_2$ in $\Dp$ such that $\th=\b_1+\b_2$, at most 
one of $\b_1, \b_2$ lies in $N(y_\a)$; moreover if $\b\perp\th$, 
then $\b\not\in N(y_\a)$. 
\par 
If $\a$ is simple, exactly one of $\b_1, \b_2$ lies in $N(y_\a)$; in particular  
$$g-2={1\over 4}(|\D|-|\D_\th|-2).$$
Conversely, assume that $y\in W$ is such that for any pair $\b_1, \b_2$ in 
$\Dp$ such that $\th=\b_1+\b_2$, exactly one of $\b_1, \b_2$ lies in $N(y)$, 
and, moreover, $\th\not\in N(y)$.
Then there exists $\a\in \Pi_\ell$ such that $y(\a)=\th$.
\endproclaim

\demo{Proof}
We have $\th\notin N(y_\a)$, therefore, by 2.3 (3),  for any pair
$\b_1, \b_2$ in $\Dp$ such that $\th=\b_1+\b_2$, at most one of $\b_1,
\b_2$ lies in $N(y_\a)$. Moreover, by 7.1 (2) and 2.3 (5), no root
orthogonal to $\th$ may belong to $N(y_\a)$. If $\a\in \Pi$,  since
$y_\a(\a)>0$, we have $N(y_\a s_\a)=N(y_\a)\cup\{\th\}$; therefore at
least one of $\b_1, \b_2$ lies in $N(y_\a s_\a)$ and hence in
$N(y_\a)$. By 7.1 (2) we have that, moreover, no root in $\D_\th$
belongs to $N(y_\a)$, therefore, using 7.1 (3), we obtain
$g-2=\ell(y_\a)={1\over 2}(|\Dp\meno \Dp_\th|-1)={1\over
4}(|\D|-|\D_\th|-2)$.
\par
Now let $y\in W$ be such that for any pair $\b_1, \b_2$ in $\Dp$ such
that $\th=\b_1+\b_2$, exactly one of $\b_1, \b_2$ lies in $N(y)$, and,
moreover, $\th\not\in N(y)$. Then $N(y)\cup \{\th\}$ still verifies
conditions ({\it a}) and ({\it b)} of 2.3 (4), therefore there exists
$y'\in W$ such that $N(y')=N(y)\cup \{\th\}$. By 2.6, $y'=y x$ with
$\ell(y')=\ell(y)+\ell(x)$, therefore $\ell(x)=1$; equivalently,
$x=s_\a$ for some $\a\in \Pi$. Then by 2.4 we obtain that $\th=y(\a)$,
hence that $\a\in\Pi_\ell$.
\endemo
\bigskip

\ni{\bf Notation.}
For $\a\in \Dp_\ell$, we set $\wa=s_0 y_\a.$ Recall  that $y_\a
s_\a=s_\th y_\a$ and $t_{\thvee}=s_0s_\th$, therefore we have
$$\wa=s_0 y_\a=t_{\thvee}y_\a s_\a.$$
\smallskip

\proclaim{Proposition 7.3} Let $\a\in \Dp_\ell$.
\roster
\item 
We have $\II(\thvee, \a)=\{\goi(\wa)\}$.
\item
For any pair $\b_1, \b_2$ in $\Dp$ such that $\th=\b_1+\b_2$, at most
one of $-\b_1+\d, -\b_2+\d$ lies in $N(\wa)$, and if $\a\in \Pi$
exactly one.
Moreover, $\a_0\in N(\wa)$, and $\b\not\in N(\wa)$ for all
$\b\in\Dap\cap\a_0^\perp$.
\item 
$N(\wa^{-1})=\{\b\in\Dp\mid (\b,\a)=-1\}\cup\{-\a+\d\}$.
\endroster
\endproclaim

\demo{Proof}\par
(1) follows from 6.1 (3) and 7.1 (2).
\par
(2) If $\th=\b_1+\b_2$ with $\b_1, \b_2\in \Dp$, then $(\b_i,\thvee)=1$ for
$i=1, 2$, so that $s_0(\b_1)=-\b_2+\d$ and $s_0(\b_2)=-\b_1+\d$. 
Since $N(\wa)=N(s_0)+s_0N(y_\a)$, we get the claim by Proposition 7.2. 
\par
(3) follows from  7.1 (5), since $N(\wa^{-1})=N(y_\a^{-1})\cup y_\a^{-1} N(s_0)$.
\endemo
\bigskip

\ni{\bf Remark.}
It is clear that for all $\t\in X^+$ we have that $\II(\t)=
\bigcup\limits_{\a\in\Dp_\ell}\II(\t, \a)$. Therefore by 7.3 (1) we
obtain that $|\II(\thvee)|=|\Dp_\ell|$. By \cite{CMP, Theorem 2.4} we obtain
that $|\Dp_\ell|=|W|/|W(\Pia\meno \Pi'_\th)|$.
\bigskip

We can apply Proposition 7.1, in particular, to the irreducible
parabolic subsystems~$\D(\tau),\,\t\in X^+$. We set 
$$\align
&\Pi_\ell(\t)=\Pi_\ell\cap\Pi(\t)\\
&\Dp_\ell(\t)=\Dp_\ell\cap\D(\t).\endalign$$
For $\a\in \Pi_\ell(\t)$ we denote by 
$y_{\t,\a}$ be the unique element of minimal length in $W(\D(\tau))$
such that $y_{\t,\a}(\a)=\th(\t)$.
\bigskip

\proclaim{Proposition 7.4}
Let $\t\in X^+$ and $\a\in \Dp_\ell(\t)$. Then $t_\t v_\t y_{\t,\a}\in
\wab(\t,\a)$. In particular, $\II(\t,\a)$ is not empty if and only if
$\a\in\Dp_\ell(\t)$.
\endproclaim

\demo{Proof}  
The \lq\lq only if \rq\rq\ part was stated in 6.10.  Conversely, if
$\a\in \Dp_\ell(\t)$, then $y_{\t, \a}(\a)=\th(\t)$ and $L_{y_{\t, \a}}
\sbsq \Pi'(\t)$. Moreover, since $v_\t y_{\t, \a}(\a)=-\th$, we have
$\th\in N(v_\t y_{\t,\a})$. By Lemma 6.5 and 6.1 (3), it
follows that  $t_\t v_\t y_{\t, \a}\in \II(\t, \a)$.
\endemo
\bigskip

\ni{\bf Remark.}
The element $y_\a$ defined in Proposition 7.1 equals $y_{\thvee, \a}$,
in the above notation. Comparing Propositions 7.3 (1) 
and 7.4 we obtain that, indeed, $v_{\thvee} y_\a=y_\a s_\a$.
\bigskip

\proclaim{Proposition 7.5}
Let $\t\in X^+$, $\a\in\Dp_\ell(\t)$, and $y\in W$ be such that 
$y(\a)=\th(\t)$ and $L_y\sbsq \Pi'(\t)$. 
Then $y_{\t,\a}\leq y$ and $y_{\t, \a}^{-1} y\in W(\Pi\cap\a^\perp)$. 
In particular 
$$\goi(t_\t v_\t y_{\t, \a})=\min \II(\t, \a).$$
\endproclaim

\demo{Proof}
By 2.8, there exists $y'\in W(\D(\t))$ such that $N(y')=N(y)\cap
\D(\t)$ and $y=y'x$ with $\ell(y)=\ell(y')+\ell(x)$. Set
$\a'={y'}^{-1}(\th(\t))$. Clearly $\a'\in \Dp(\t)$, and since
$L_{y'}\sbsq \Pi'(\t)$, by 7.1 (2) we obtain that $y'=y_{\t, \a'}$.
Then  $x=y_{\t, \a'}^{-1} y$; hence $y=y_{\t, \a'} x$ and
$N(y)=N(y_{\t,\a'}) \cup y_{\t, \a'} N(x)$.
We prove that $x$ is a product of reflections with
respect to simple roots orthogonal to $\a'$. Let $x=x_1\cdots x_k$ with
$x_i$ simple reflection for $1\leq i \leq k$. Assume by contradiction
that $x_j=s_\b$ with $\b\not\perp \a'$ for some $j\in \{1, \dots, k\}$,
and $x_i(\a')= \a'$ for $i<j$. Set $x'= x_1\cdots x_{j-1}$, and 
$y'=y_{\t, \a'} x'$. Then $y'(\b)\in N(y)\meno N(y_{\t, \a'})$, 
hence $y'(\b)\notin \D(\t)$; in particular $y'(\b)\not \leq \th(\t)$. We
have
$y'(\a')=\th(\t)$ and, since $(\b,\a')\ne 0$, $(y'(\b),\th(\t))\ne 0$,
hence
either $y'(\b)-\th(\t) $, or $y'(\b)+\th(\t) $ is a root. In the second
case,
by 6.8 (4) we obtain that  $n_\g(y'(\b))=0$ for all $\g\in \Pi'(\t)$:
this is impossible, since $y'(\b)\in N(y)$ and $L_y\sbsq \Pi'(\t)$. 
In the first case, since $y'(\b)\not \leq \th(\t)$, we obtain that 
$y'(\b)-\th(\t) >0$; moreover, since $y'(\b)\in N(y)$ and $\th(\t)\notin
N(y)$, 
by 2.3 (3) 
$y'(\b)-\th(\t)\in N(y)$. As above this leads to a contradiction, since   
$n_\g(y'(\b)-\th(\t))=0$ for all $\g\in \Pi'(\t)$.
Thus we have in particular that $x(\a')=\a'$, hence $y(\a')=\th(\t)$; it
follows that $\a'=\a$. By Proposition 6.6 we obtain that 
$\goi(t_\t v_\t y_{\t, \a})=\min \II(\t, \a)$.
\endemo
\bigskip

We set, for $\t\in X^+$, $\a\in\Dp_\ell(\t)$ $$\alignat2
&W_\a &&=W(\Pi\cap \a^\perp),\\
&P_{\t,\a} &&=\{\b\in \Pi\meno \Pi(\t)\mid \b\not\perp \D(\t), \
\b\perp \a\}.\endalignat$$
\smallskip

\proclaim{Lemma 7.6}
Let $\t\in X^+$, $\a\in\Dp_\ell(\t)$, and $y\in W$ be such that 
$y(\a)=\th(\t)$ and $L_y\sbsq \Pi'(\t)$. 
Set $x=y_{\t,\a}^{-1} y$.
Then  $x\in W_\a$ and $L_x\sbsq P_{\t, \a}$. 
\par
Conversely, if $x\in W_\a$ and $L_x\sbsq P_{\t,\a}$,
then  $y_{\t,\a}x(\a)=\th(\t)$ and $L_{y_{\t,\a}x}\sbsq \Pi'(\t)$.
Moreover, $\ell(y_{\t,\a}x)=\ell(y_{\t,\a})+\ell(x)$.
\endproclaim

\demo{Proof}
By Proposition 7.5, we have that $x\in W_\a$, hence $N(x)\sbsq
\D(\Pi\cap\a^\perp)$; moreover $N(y)=N(y_{\t,\a})\cup y_{\t,\a} N(x)$.
Therefore it suffices to prove that for all $\b\in \Pi\cap\a^\perp
\meno P_{\t, \a}$ we have that $\b\not\in N(x)$. Assume $\b\in
\Pi\cap\a^\perp \meno P_{\t, \a}$; this means that either $\b\in \Pi(\t)$, 
or $\b\not\in \Pi(\t)$ and $\b\perp \D(\t)$. First assume $\b\not\in
\Pi(\t)$ and $\b\perp\D(\t)$. Then, since $y_{\t,\a}\in W(\D(\t))$, we
obtain that $y_{\t,\a}(\b)=\b$. By assumption, $\b\not\in N(y)$, hence
$\b\not\in N(x)$. Next assume $\b\in \Pi(\t)$. Then $y_{\t,\a}(\b) \in
\D(\t)$ and, since $\b\perp\a$, $y_{\t,\a}(\b)\perp\th(\t)$. It
follows that $\supp(y_{\t,\a}(\b)) \cap\Pi'(\t)=\emptyset$; therefore
$y_{\t,\a}(\b)\not\in N(y)$, whence $\b\not\in N(x)$.
\par
Conversely, let  $x\in W_\a$ be such that $L_x\sbsq P_{\t, \a}$. It is
clear that $y_{\t,\a}x(\a)=\th(\t)$.  We use 2.3 (5): since
$N(y_{\t,\a})\subseteq \Pi'(\t)$, it suffices to prove that, for all
$\b\in N(x)$, $\supp(y_{\t,\a}(\b))\cap \Pi'(\t)\ne\emptyset$, and
$y_{\t,\a}(\b)>0$.  By assumption, for each $\b\in N(x)$,
$\supp(\b)\cap P_{\t,\a}\ne \emptyset$. Since $y_{\t,\a}\in W(\D(\t))$,
$y_{\t,\a}$ acts on any root by summing linear combinations of roots
in $\Pi(\t)$. Since $P(\t,\a)\cap \Pi(\t)=\emptyset$, we obtain that
$\supp(y_{\t,\a}(\b))\cap P_{\t,\a}\ne \emptyset$; moreover
$y_{\t,\a}(\b)>0$. But $\b\perp \a$, hence
$y_{\t,\a}(\b)\perp\th(\t)$.  Since $\th(\t)$ has negative scalar
product with all roots in $P_{\t,\a}$, we obtain that
$\supp(y_{\t,\a}(\b))$ also contains some root having positive scalar
product with $\th(\t)$. By Lemma 5.5 (3) it follows that
$\supp(y_{\t,\a}(\b))\cap \Pi'(\t)\ne\emptyset$. This implies the
thesis.
\endemo
\bigskip

Notice that if $\b\in \Pi\meno \Pi(\t)$ and $\b\not\perp\D(\t)$, then
$(\b,\th(\t))<0$. Let $\{\D=\D_1\sps\cdots \sps\D_k\}$ be the good
chain of $\t$. By Lemma 5.5, for all $\b\in \Pi$ we have $(\b,
\th(\t))<0$ if and only if $\b\in \Pi'_{k-1}$, where
$\Pi'_{k-1}=\{\b\in\Pi\mid(\b,\th_{k-1})\ne 0\}$, $\th_{k-1}$ being the
highest root of $\Pi_{k-1}$. Thus we obtain that $P_{\t,
\a}=\Pi'_{k-1}\cap\a^\perp$.
We set, for $\a\in\Dp_\ell$ and $\t\in X^+$
$$W_{\a, \t}=W((\Pi\cap \a^\perp)\meno \Pi'_{k-1}).$$ 
As seen in 2.9, $W_{\a, \t}\backslash W_\a $ is a poset with minimum
and maximum. If  $w_0^\a$ is the longest element of $W_\a$, and
$w_0^{\a,\t}$ is the longest element of  $W_{\a, \t}$, then
$w_0^{\a,\t}w_0^\a$ corresponds to the maximal element of $W_{\a, \t}\backslash
W_\a$. We set 
$$\check y_{\a,\t}=y_{\a,\t} w_0^{\a,\t}w_0^\a.$$
We can now describe the poset structure of $\II(\t, \a)$.
\bigskip

\proclaim{Proposition 7.7}
Let $\t\in X^+$,  and $\a\in \Dp_\ell(\t)$.
Then the poset $\II(\t,\a)$ is
isomorphic to $W_{\a, \t}\backslash W_\a$.
In particular $\II(\t,\a)$ has a unique maximal element and a unique
minimal element. Moreover,  
$$\max \II(\t,\a)=\goi(t_\t v_\t \check y_{\a,\t}).$$
\endproclaim

\demo{Proof}
Assume $\goi, \goj\in \II(\t,\a)$. Then by Proposition 6.6, $w_\goi=t_\t
v_\t y$, 
$w_\goj=t_\t v_\t y'$, with $L_y, L_{y'}\sbsq \Pi'(\t)$. Moreover
$\goi\sbsq \goj$
if and only if $y\leq y'$. By Lemma 7.6, there exist $x, x'\in W_\a$ such
that
$L_x, L_{x'}\sbsq P_{\t, \a}$, $y=y_{\t, \a}x$, and $y'=y_{\t, \a}x'$.
Moreover, 
$y\leq y'$ if and only if $x\leq x'$. By 2.9, this implies the claim.
\endemo
\vskip 1cm

\heading \S8 Maximal elements in $\II(\t,\a)$\endheading
\medskip
In the previous section we have shown that $\II(\t,\a)$ is not empty
if and only if $\a\in\Dp_\ell(\t)$ and in this case it has unique minimal and
maximal elements, but at this point we have no information
about the relationships between $\II(\t,\a)$ and $\II(\t',\a)$. In
this section we shall prove that $\max\II(\t,\a) \subseteq
\max\II(\t',\a)$ if and only if $\tau\leq\t'$. We remark that if
$\t\leq \t'$, then $\max \II(\t,\a)$ and $\min\II(\t',\a)$ may be not
comparable.
\medskip

Let $\a\in \Dp_\ell$. We set 
$$\adj(\a)=\{\a'\in \Pi\mid (\a, \a')\ne 0\};$$ 
$$\Phi_*(\a)=\{\b\in \D\mid n_{\a'}(\b)=m_{\a'}, \text{\ for
all\ }\a'\in\adj(\a)\};$$ 
$$\goi_*(\a)=\sum\limits_{\b\in\Phi_*(\a)}\gog_\b.$$
Notice that $\adj(\a)\cap\supp(\a)\ne 0$, since $(\a,\a)\ne 0$.
\par
Consider the parabolic subsystems of $\Da$, $\D$, generated by all
simple roots which are orthogonal to $\a$,  $\D(\Pia\cap \a^\perp)$
and $\D(\Pi\cap\a^\perp)$, respectively. Recall that, for $\a\in
\Dp_\ell$, we have set $\wp=W(\Pi\cap\a^\perp)$. We also set:
$$\wap=W(\widehat\Pi\cap\a\perp),\quad \dap=\D(\widehat\Pi\cap\a\perp),
\quad \dperp=\D(\Pi\cap\a\perp).$$
It is clear that $\dperp\sbsq \dap$ and that equality holds if and
only if $\a\not\perp\th$.
\smallskip

\proclaim{Lemma 8.1} {\rm (1)} 
For all $\a\in \Dp_\ell$, $\goi_*(\a)$ is an abelian ideal of $\gob$.
\parni
{\rm (2)}
Assume  $\a\perp\th$, and set $w_{*,\a}=w_{\goi_*(\a)}$. Then  
$$N(w_{*,\a})=-\Phi_*(\a)+\d=\dap^+\meno
\dperp^+;$$ in particular 
$w_{*,\a}\in \wab\cap \wap$.
\parni 
{\rm (3)}
If $x\in \wap$ and $L_x\sbsq \{\a_0\}$, 
then $x\in \wab$ and $x\leq w_{*,\a}$. In particular, 
for all $x\in \wab\cap \wap$, we have $x\leq w_{*,\a}$. 
\endproclaim

\demo{Proof}
\parni
(1) It is clear that $\Phi_*(\a)\sbsq\Dp$ and is an abelian set, i.e., 
for all $\b, \b'\in \Phi_*(\a)$, $\b+\b'$ is not a root. It also clear 
that if $\b\in\Phi_*(\a)$ and $\b'\in \Dp$ is such that $\b'\geq \b$, 
then $\b'\in\Phi_*(\a)$. This proves that $\goi_*(\a)$ is an abelian 
ideal of $\gob$.
\parni
(2) The first equality follows from section 3 and holds even if 
$\a\not\perp\th$. We prove the second equality. 
\par
Assume that $\a_0\perp \a$ and $\b\in \dap^+\meno \dperp^+$. 
Then $\b=c\a_0+\sum\limits_{\g\in\Pi\cap\a^\perp}c_\g \g$, with $c,
c_\g\in \nat$ and $c>0$, hence $\b=\ov \b+c\d$ with   
$\ov\b=-c\th+\sum\limits_{\g\in\Pi\cap\a^\perp}c_\g \g$. Now  
$n_{\a'}(\ov\b)= -c m_{\a'}$ for all $\a'\in \adj(\a)$.
It follows that  $c=1$, hence that $-\ov\b\in \Phi_*(\a)$, and
$\b\in -\Phi_*(\a)+\d$.
\par
Conversely, if $\ov\b\in \Dp$ is such that $n_{\a'}=m_\a$ for all
$\a'\in \adj(\a)$, then $-\b+\d-\a_0$ is either $0$, or a non negative
linear combination of roots in $\Pi\cap\a^\perp$, hence $-\b+\d\in
\dap^+ \meno \dperp^+$.
\parni
(3) Assume $x\in \wap$ and $L_x\sbsq \{\a_0\}$. Then $N(x)\sbsq
\dap^+$; moreover, by 2.3 (3) we have that $N(x)\cap\D=\emptyset$,
hence $N(x)\sbsq\dap^+\meno\dperp^+$. By (2) and 2.6 we obtain
that $x\leq w_{*,\a}$, and hence that $x\in \wab$.
\endemo
\bigskip

\ni{\bf Remark.}
Since for all $w\in \wab$ we have $L_w\sbsq\{\a_0\}$, by Lemma 
8.1 (3) and 2.9 we obtain that the poset $\wab\cap\wap$ is isomorphic
to $\wp\backslash \wap$.
\bigskip 

\proclaim{Lemma 8.2}
Let  $\{\D_1=\D\sps\cdots\sps\D_k\}$ be a good chain, 
and $\Pi_i=\D_i\cap \Pi$ for $i=1, \dots, k$. 
Assume that $k\geq 2$ and $\a\in \Pi_k$. Then $n_\a(\th_{k-1})=m_\a$.
\endproclaim

\demo{Proof}
We assume by contradiction that  $n_\a(\th_{k-1})<m_\a$. 
We shall prove that then for all $\b\in \Pi_k$ we 
have that $n_\b(\th_{k-1})<m_\b$: this is impossible since, 
by 6.8 (4), for $\b\in \Pi'_k$ we have $m_\b=n_\b(\th_k)\leq
n_\b(\th_{k-1})$. 
\par
Assume $\b\in \Pi_k$. We first prove that, if $n_{\b'}(\th_{k-1})<m_{\b'}$
for some $\b'\in\adj(\b)$, then $n_\b(\th_{k-1})<m_\b$. 
Since $\b\in\D_k$, we have $(\b,\th_i)=0$ for all $i\in [k-1]$, hence 
$0=(\b,\th)=\sum\limits_{\b'\in\adj(\b)} m_{\b'}(\b, \b')$, and
$0=(\b,\th_{k-1})=\sum\limits_{\b'\in\adj(\b)} n_{\b'}(\th_{k-1})(\b, \b')$. 
Since $(\b, \b)>0$ and $(\b, \b')<0$ for all $\b'\ne \b$ in $\adj(\b)$, 
we obtain that, if $n_{\b'}(\th_{k-1})<m_{\b'}$ for some  $\b'\ne \b$ in $\adj(\b)$, 
then $n_\b<m_\b$.
\par
Now, for all $\b\in\Pi_k$, there exists a sequence, $\b_1=\a, \dots,
\b_h=\b$, in $\Pi_k$, such that $\b_{i+1}\in\adj(\b_{i})$ for $i=1,\dots,
h-1$. If $n_\a(\th_{k-1})<m_\a$, we obtain inductively that
$n_\b(\th_{k-1})<m_\b$.
\endemo
\bigskip

Let $\a\in \Dp_\ell$. 
We recall that, by Remark 5.7 (3), there exists a unique maximal good
chain $\{\D=\D_1\sps \cdots\sps \D_{k}\}$ such that $\a\in \D_i$ for
$1\leq i\leq k$: we shall call this chain the {\it maximal good chain
of $\a$}. Any good chain the elements of which contain $\a$ is an
initial section of the maximal one.
\bigskip

\ni
{\bf Notation.} Let $\a\in\Dp_\ell$. We denote by $\{\D=\D_1(\a)\sps
\cdots\sps \D_{k(\a)}(\a)\}$ the maximal good chain of $\a$, and  for
$1\leq h\leq k(\a)$, we denote by $\th_h(\a)$ the highest root of
$\D_h(\a)$.  Moreover we set $\t_0(\a)=0$, and  
$$\t_h(\a)=\sum\limits_{i=1}^h\th_i(\a)^\vee,\quad 
\text{for $1\leq h\leq k(\a)$}, \quad X_\a=\{\t_h(\a)\mid
1\leq h\leq k(\a)\}.$$ 
It is clear that if $\a\in \Dp_\ell$, then $\t\in X_\a$ if 
and only if $\a\in \Dp_\ell(\t)$.
Thus by Remark 5.7 and Proposition 7.4 we have 
that $\II(\t,\a)\ne\emptyset$ if and only if $\t\in X_\a$.
\bigskip

\proclaim{Proposition 8.3}
Let $\a\in \Dp_\ell$ and assume that $k(\a)\geq 2$. 
Then $-\th_h(\a)+\d\in\dap$, for all $h\in \{1, \dots, k(\a)-1\}$. 
Equivalently, for $1\leq h<k(\a)$, we have $n_\b(\th_h(\a))=m_\b$, for all 
$\b\in\adj(\a)$.
\endproclaim

\demo{Proof}
Set $k=k(\a)$. We shall prove that for all $\b\in\adj(\a)$ we have
$n_\b(\th_{k-1})=m_\b$; this clearly implies that $n_\b(\th_h)=m_\b$
for all $h\in \{1, \dots, k-1\}$, which is equivalent to the thesis.
Let $\b\in \adj(\a)$. If $\b\in \Pi_k$, then the claim follows directly 
from Lemma 8.2. Assume $\b\not\in \Pi_k$. We first show that $\b\in  \Pi_{k-1}$. 
If not, consider the maximal $i\in \{1,\dots, k-2\}$ such that 
$\b\in \Pi_i$. Since $\a\in \Pi_{i+1}$, $\b$ cannot 
belong to any irreducible component of $\D_i\cap\th_i^\perp$ different from 
$\D_{i+1}$; it follows that $\b\in\Pi'_i$. But then, by property 
(2) of Definition 5.4, $\b\perp\D_{i+2}$: this is impossible, 
since $\a\in\D_{i+2}$ and $\b\in\adj(\a)$. Thus we have $\b\in
\Pi_{k-1}\meno\Pi_k$ and $\b\not\perp\D_k$: this implies that
$\b\in \Pi'_{k-1}$, hence, by 6.8 (4), $n_\b(\th_{k-1})=m_\b$.
\endemo
\bigskip

In the next Lemma we consider the element $\wa$ defined in section 7. 
 
\proclaim{Lemma 8.4}
$(1)$ Let $\a\in\Dp_\ell$ and $x\in \wab\cap \wap$. Then $\wa x\in \wab$
and $\goi (\wa x)\spsq \goi(\wa)$. Moreover, $\dim\goi(\wa x)=
\ell(\wa)+\ell(x)$.
\parni
$(2)$ Let $\a\in\Dp_\ell$ and $x\in \hatw_\a$ be such that $\wa x\in
\wab$. Then $x\in \wab\cap \hatw_\a$.
\endproclaim

\demo{Proof}
(1) If $x\in \wab\cap \wap$, then we have that $N(x)\sbsq
\Dm+\d$, and $N(x)\perp \a$.  Consider $z\in
N(x),\,z=-\b+\d,\,\b\in\Dp$. We have $\b\perp\a$, so, by 7.1 (5),
$\b\notin N(y_\a^{-1})$; it follows  that $y_\a N(x)\sbsq \Dm+\d$.
Moreover, $y_\a N(x)\perp \th$, therefore $\wa N(x) =s_0 y_\a N(x)=
y_\a N(x)$. In particular, $\wa N(x)\in \Dap$, hence by  2.7 $\wa \leq
\wa x$. Moreover, $N(\wa x)=N(\wa)\cup \wa N(x)\sbsq \Dm+\d$, hence
$\wa x\in \wab$. By 3.2 we obtain that $\goi(\wa)\sbsq \goi(\wa x)$.
Moreover, we have $\dim\goi(\wa x)=\ell(\wa x)=\ell(\wa)+\ell(x)$. 
\parni
(2) By 8.1 (3) it suffices to prove that $L_x\sbsq \{\a_0\}$. If
$L_x\not\sbsq\{\a_0\}$, then there exists $\b\in\Pi\cap \a^\perp$ such
that $\b\in N(x)$. By Proposition 7.3 (3) we have $\wa(\b)>0$, hence
$\wa(\b)\in N(\wa x)\sbsq \Dm+\d$. It follows that $\wa s_\b\in \wab$
and $\goi(\wa s_\b)\supsetneq \goi(\wa)$. But $\goi(\wa s_\b)\in
\II(\thvee, \a)$, and, by 7.3 (1), $\II(\thvee,\a)=\{\i(\wa)\}$, thus
we get a contradiction.
\endemo
\bigskip

\proclaim{Lemma 8.5}
Let $\a\in \Dp_\ell$. Assume $\t\in X$,  $x\in \wab(\t)\cap
\wap$ and set $2h=\frac{(\th,\th)}{2}(\t,\t)$. Then:
\roster 
\item 
$h\in \{0, \dots, k(\a)-1\}$ and $\wa x\in\wab(\t_{h+1}(\a), \a)$;
\item 
$y_\a(\t)=\th_2(\a)^\vee+\cdots+\th_{h+1}(\a)^\vee$.
\endroster
\endproclaim

\demo{Proof}
By 8.4 (1), $\wa x\in \wab$; moreover, it is clear that 
$\wa x(\a)=-\th+2\d$. Let $x=t_\t v$, with $v\in W$. 
Then $\wa x\in \II(\thvee+y_\a(\t), \a)$. By Proposition 7.4 and
Remark 5.7 (3), we obtain that  $\thvee+y_\a(\t)=\t_j(\a)$ for some
$j\in \{1,\dots, k(\a)\}$. Since $\a\perp\t$, we have that $\th\perp
y_\a(\t)$, hence ${(\th, \th)\over
2}(\thvee+y_\a(\t),\thvee+y_\a(\t))=2(h+1)$. By Remark 5.7 (1), it
follows that $\thvee+y_\a(\t)=\t_{h+1}(\a)$ and hence that $1\leq
h+1\leq k(\a)$, so that  $h\in \{0, \dots, k(\a)-1\}$. Moreover
$y_\a(\t)=\th_2(\a)^\vee+\cdots+\th_{h+1}(\a)^\vee$.
\endemo
\bigskip

\ni
{\bf Notation.} 
Let $\a\in\Dp_\ell$. If $k(\a)>1$ and $h\in[k(\a)-1]$, we set
$$\Psi_h(\a)=\{\b\in \D_h(\a)\mid \b=\th_h(\a), \text{ or \ }
\th_h(\a)-\b\in \dperp\};$$
$$\Phi_h(\a)=\bigcup\limits_{i=1}^h\Psi_i(\a);\quad
\goi_h(\a)=\sum\limits_{\b\in\Phi_h(\a)}\gog_\b.$$
Moreover, we set $\Phi_0(\a)=\emptyset$, and $\goi_0(\a)$ to be the
zero ideal. By Proposition 8.3 we directly obtain that
$\Phi_h(\a)\sbsq\Phi_*(\a)$ for all $h\in[k(\a)-1]$. We remark that
the subsets $\Psi_h(\a)$, $h\in[k(\a)-1]$, are pairwise disjoint.
\bigskip

\proclaim{Lemma 8.6}
Let $\a\in \Dp_\ell$. Then $\goi_h(\a)$ is an abelian ideal of $\gob$,
for $0\leq h <k(\a)$. Set $\wah=w_{\goi_h(\a)}$. Then
$\wah\in\wab(\t_h(\a))\cap \wap$. In particular
$\wa\wah\in \wab(\t_{h+1}(\a),\a)$ and
$\dim\goi(\wa\wah)=\ell(\wa)+\ell(\wah)$.
\endproclaim

\demo{Proof}
We set $k=k(\a)$; $\D_h=\D_h(\a)$, $\th_h=\th_h(\a)$, $\t_h=\t_h(\a)$,  
and $\Phi_h=\Phi_h(\a)$, for $0\leq h <k$.
\par
If $h=0$, the claim is trivial. Assume $k\geq 2$ and $1\leq h\leq
k-1$. Since $\Phi_h(\a)\sbsq \Phi_*(\a)$, by Lemma 8.1 (1) we obtain
that $\goi_h(\a)$ is an abelian subalgebra.
\par
Then let $\b\in \Phi_h$ and  $\b'\in \Dp$ be such that and
$\b'\geq\b$. We shall prove that $\b'\in \Phi_h$; this will imply that
$\goi_h(\a)$ is an abelian ideal.
It is clear that $n_\g(\b')=m_\g$ for all $\g\in\adj(\a)$. Let $i\in
[h]$ be such that $\b\in\D_i$ and $\th_i-\b \in \dperp\cup\{0\}$.
Assume first $\b'\in \D_i$. Then, since $\th_i$ is dominant in $\D_i$,
we have $(\th_i,\b')\geq (\th_i,\b)>0$. Hence either $\b'=\th_i$, or
$\th_i-\b'\in\D$. In the latter case, since $\th_i-\b'\leq \th_i-\b$,
we have that $\th_i-\b' \in \dperp$; therefore $\b'\in \Phi_h$.
Next assume $\b'\notin\D_i$ and let $j$, $1\leq j\leq i-1$,  be such
that $\b'\in \D_j$ and $\b'\not\in \D_{j+1}$. If  $\b'\perp \th_j$,
then $\b'$ is included in some connected component of
$\D_j\cap\th_j^\perp$. But, for all $\a'\in\adj(\a)$, we have
$n_{\a'}(\b')=m_{\a'}>0$; moreover $\adj(\a)\cap\supp(\a)\ne
\emptyset$, and $\supp(\a)\in \D_{j+1}$, hence we obtain that
$\b'\in\D_{j+1}$, a contradiction.
Hence $\b'\not\perp\th_j$. It follows that either $\b'=\th_j$, or
$\th_j-\b'\in \D$, and, indeed, $\th_j-\b'\in \dperp$. Thus $\b'\in
\Phi_h$. This concludes the proof that $\goi_h(\a)\in \II$.
\par
By definition, $\th_i\in \Phi_h$ for all $i\in [h]$; moreover
$\b\not\in\Phi_h$ for all $\b\in\D_h\cap\th_h^\perp$. It follows that
$-\th_i+\d\in N(\wah)$ for all $i\in [h]$, and $-\b+\d\not\in N(\wah)$
for all $\b\in\D_h\cap\th_h^\perp$. By Proposition 6.2, it follows
that $\wah \in \wab(\t_h)$.
\par
If $\b\in \Psi_j(\a)$, for some $j\in [h]$, then, by definition,
$-\b+\d-(-\th_j+\d)\in \dperp$, and since, by 8.3, $-\th_j+\d\in
\dap$, we obtain that $-\b+\d\in \dap$. Since $N(\wah)=
\{-\b+\d\mid\b\in \Phi_h\}$, we obtain that $N(\wah)\sbsq \dap$, hence
that $\wah\in \wap$.
\par
By  8.5 and 8.4 (1) we finally obtain that $\wa\wah\in
\wab(\t_{h+1}(\a),\a)$, and\break $\dim\goi(\wa\wah)=\ell(\wa)+
\ell(\wah)$.
\endemo
\bigskip

\ni
{\bf Remark.} For $\a\in\Dp_\ell$ and $0\leq h<k(\a)$, let 
$\wah=w_{\goi_h(\a)}$, as in the statement of 8.6. 
Notice that $\goi_0(\a)$ is the zero ideal, hence $\check w_{0,\a}=1$.
Notice, moreover, that, for $h>0$
$$N(\wah)=-\Phi_h(\a)+\d.$$
\bigskip

\proclaim{Proposition 8.7}
Let $\a\in\Dp_\ell$, $0\leq h< k(\a)$, and $\goi\in \II(\t_{h+1}(\a), \a)$.
Then $\goi\supseteq \goi(\wa)$, and $\wa^{-1}w_\goi\in \wab(\t_h(\a))
\cap \hatw_\a$.
\endproclaim

\demo{Proof} 
Set $\t=\t_{h+1}(\a)$. By Lemma 8.6 we have that $\wa \wah\in \wab(\t,
\a)$. Set $\goj=\min\II(\t,\a)$. By Proposition 7.5, for each $w\in
\wab(\t, \a)$ we have that $w=w_\goj y$ with $y\in \wp$. In
particular, for all $w, w'\in \wab(\t,\a)$ we have that $w^{-1}w'\in
\wp$. Thus $w_\goi=\wa \wah y $, with $y\in \wp$, and
therefore $\wah y\in \wap$. Then by Lemma 8.4 (2) we obtain that $\wah
y\in \wab\cap \hatw_\a$, and it is clear that $\wah y\in \wab(\t_h(\a),
\a)$. Moreover, by 8.4 (1), $\goi\supseteq \goi(\wa)$.
\endemo
\bigskip

For $\a\in \Dp_\ell$ we set 
$$\II(\a)=\bigcup\limits_{\t\in X_\a}\II(\t, \a).$$
These posets have been studied also in \cite{P}, \cite{S}. If $\goi\in\II$, then $\goi\in\II(\a)$ if and only if
$w_\goi(\a)=-\th+2\d$.
We can now obtain the following results.
\bigskip

\proclaim{Corollary 8.8} (Suter)
The poset $\II(\a)$ is isomorphic to $W_\a\backslash\hatw_\a$. In particular, 
it has a unique minimal element and a unique maximal element. 
Moreover, 
$$\min\II(\a)=\goi(\wa),\qquad \max\II(\a)=\goi(\wa w_{*, \a}).$$
\parni
\endproclaim

\demo{Proof}
By Proposition 8.7 and Lemma 8.4, $\goi\in \II(\a)$ if and only if
$w_\goi=\wa x$ with $x\in \wab\cap\wap$; moreover, in this case,
$N(\wa x)=N(\wa)\cup \wa N(x)$. Therefore, if $\goi, \goj\in \II(\a)$, 
$w_\goi=\wa x$, and $w_\goj=\wa y$, with $x, y\in \wab\cap\wap$, 
then $\goi\sbsq\goj$ if and only if $x\leq y$. By Lemma 8.1
and the subsequent remark we get the claim.
\endemo
\bigskip

\proclaim{Corollary 8.9} 
Assume $\goi\in \II,\,\a\in\Dp_\ell$ and $\t\in X_\a$. 
\roster 
\item 
$\goi\in \II(\a)$ if and only if $\min\II(\a)\sbsq \goi\sbsq \max\II(\a)$.
\item 
$\goi\in \II(\t,\a)$ if and only if $\min\II(\t, \a)\sbsq \goi\sbsq
\max\II(\t, \a)$.
\endroster
\endproclaim

\demo{Proof} 
(1) The \lq\lq only if\rq\rq\ part is a trivial consequence of
Corollary 8.8. 
Conversely, assume that  $\min\II(\a)\sbsq \goi\sbsq \max\II(\a)$. Then
there exists $x\in W$ such that $w_\goi= \wa x$ and
$N(w_\goi)=N(\wa)\cup \wa N(x)$. Since $N(w_\goi)\sbsq N(\wa
w_{*,\a})$, it follows that $N(x)\sbsq N(w_{*, \a})\sbsq \Da_\a$;
therefore $x(\a)=\a$ and $w_\goi(\a)=\wa(\a)=-\th+2\d$. 
\parni
(2) The \lq\lq only if\rq\rq\ part follows from Propositions 7.5 and
7.7. Conversely, assume $\min\II(\t, \a)\sbsq \goi\sbsq \max\II(\t,
\a)$, and let $\t=\t_h(\a)$, $h\in [k(\a)]$. By item (1), $\goi\in
\II(\a)$, hence $\goi\in \II(\t_{h'}(\a),\a)$, for some $h'\in
[k(\a)]$.  By Proposition 6.2, $\gog_{\th_i(\a)}\sbs \goi$ if
and only if  
$i\in [h']$; similarly, $\gog_{\th_i(\a)}$ is contained in 
$\min
\II(\t, \a)$ and 
$\max\II(\t, \a)$, hence in $\goi$, if and only if $i\in
[h]$.  Hence $h=h'$.
\endemo
\bigskip

We set $$\gom(\a)=\max\II(\a),$$ thus $\gom(\a)=\goi(\wa w_{*,\a})$.
\bigskip

\proclaim{Corollary 8.10} (Suter). For all $\a\in\Dp_\ell$
$$\dim\gom(\a)=g-1+{1\over
2}\left(|\dap|-|\dperp|\right).$$ 
If, moreover, $\a\perp\th$, denote by $\Da_1(\a)$ the connected component 
of $\dap$ which includes $\a_0$ and set $\Da_1^0(\a)=
\Da_1(\a)\cap\D$. Then 
$$\dim\gom(\a)=g-1+{1\over 2}\left(|\Da_1(\a)|-|\Da_1^0(\a)|\right).$$
\endproclaim
\bigskip

By Proposition 7.7, if $\II(\t, \a)$ is not empty, then it has a
unique maximal element. For all $1\leq h\leq k(\a)$, we set
$$\gom_h(\a)=\max \II(\t_h(\a), \a).$$ We can now explicitly describe
$\gom_h(\a)$.
\bigskip

\proclaim{Proposition 8.11}
For all $\a\in \Dp_\ell$, and  $0\leq h < k(\a)$, we have 
$$\gom_{h+1}(\a)=\goi(\wa \wah).$$ 
\endproclaim

\demo{Proof}
Set $\t^{-}=\t_h(\a)$. By Proposition 8.7,
$w_{\gom_{h+1}(\a)}=\wa x$, with $x\in \wab(\t^-)\cap \wap$. By
Lemma 8.6 it suffices to prove that $\wa x\leq \wa\wah$, and this is
equivalent to prove that $N(x)\sbsq N(\wah)$. Recall that
$N(\wah)=\{-\b+\d\mid \b\in \Phi_h(\a)\}$. Assume $\b\in N(x)$ and
$\b=-\ov \b+\d$, with $\ov \b\in \Dp$. We shall prove that $\ov\b\in
\Phi_h(\a)$.
\par
We set $\D_i=\D_i(\a)$, and $\th_i=\th_i(\a)$, for $1\leq i \leq h$.
Since $\D_h=\D(\t^{-})$, by  Proposition 6.2 we have that $\ov\b\not\in
\D_h\cap \th_h^\perp$, in particular $\ov\b\not\in\D_{h+1}$. Let
$i\in[h]$ be such that $\ov\b\in \D_i$ and $\ov\b\not\in \D_{i+1}$;
then $\ov\b\in\D_j$ for $1\leq j\leq i$. We first prove that
there exists $j\in \{1, \ldots, i\}$ such that $(\ov\b, \th_j)>0$. 
We assume, by contradiction, that $(\ov\b, \th_j)=0$ for $1\leq j\leq
i$. Since $\ov\b\not\in \D_{i+1}$, we obtain that $\ov\b$ belongs to a
connected component of $\D_i\cap\th_i^\perp$ distinct from $\D_{i+1}$,
in particular $\ov\b\perp\D_{i+1}$, and hence $\ov\b\perp\th_j$ also
for $j\geq i+1$. But this implies that $\ov\b\perp \t^{-}$, which is
impossible, by 6.1 (4). Therefore there exists $j\in \{1, \ldots, i\}$
such that  $(\ov\b, \th_j)>0$. Then either $\ov\b=\th_j$, or
$\th_j-\ov\b\in\Dp$. Moreover, since $-\th_j+\d$ and $\b$ belong to 
$\dap$,  we have that
$\th_j-\ov\b=\b-(-\th_j+\d)\in \dap\cup\{0\}$. But
$\th_j-\ov\b\in \D\cup\{0\}$, therefore $\th_j-\ov\b\in
\dperp\cup\{0\}$. This means that $\ov \b\in \Phi_h(\a)$.
\endemo
\bigskip

For any good chain $\{\D=\D_1\sbs\cdots \sbs \D_k\}$, we have that 
$\Pi_i\cup\{-\th_i+\d\}$ is a basis for a subsystem $\Da_i$ of $\Da$
which can be identified with the affine root system of $\D_i$.
\par
Let $\a\in \Dp_\ell$, $k=k(\a)$, and $\D_i=\D_i(\a)$, for $i\in [k]$.  
We consider the parabolic subsystem of $\Da_i$ generated by
$(\Pi_i\cup\{-\th_i+\d\}) \cap \a^\perp$: this is a finite and not
necessarily irreducible system. If $1\leq i<k$, then
$-\th_i+\d\perp\a$, and we denote by $\Da_i(\a)$ the irreducible component
of $\D((\Pi_i\cup\{-\th_i+\d\})\cap \a^\perp)$ which includes
$-\th_i+\d$. By 8.3,  $\Da_i(\a)\sbsq \dap$. 
Notice that if $\a\perp\th$, then $\Da_1(\a)$ is equal to the 
system $\Da_1$ which occurs in 8.10. 
\par
If $k(\a)>1$ and $1\leq h\leq k(\a)-1$, we set 
$$\ov\Psi_h(\a)=\{\b\in\Da_h(\a)^+\mid (\b,(-\th_h+\d)^\vee)=-1\}
\cup\{0\}.$$
\smallskip

\proclaim{Lemma 8.12}
Let $1\leq h<k(\a)$. We have 
$$\ov\Psi_h(\a)=\{\th_h-\b\mid \b\in \Psi_h(\a)\}.$$
In particular, 
$$|\Psi_h(\a)|=|\ov\Psi_h(\a)|=g_h(\a)-1, $$
where $g_h(\a)$ is the dual Coxeter number of $\Da_h(\a)$.
\endproclaim

\demo{Proof}
Since $\Da_h(\a)^+\sbsq \dap^+$, by Lemma 8.1 we
obtain that $\Da_h(\a)^+\meno \D\sbsq
-\Phi_*(\a)$ $+\d$; in particular $\Da_h(\a)^+\meno \D$ is an abelian
set.
If $\b\in\Da_h(\a)^+$ and $(\b, (-\th_h+\d)^\vee)=-1$, then
$\b+(-\th_h+\d)$ is a root, and therefore $\b\in \D$.
Then it is clear that  $\th_h-\b\in \D_h$, hence, by definition, 
$\th_h-\b\in \Psi_h(\a)$.
\par
Conversely, assume $\b\in\Psi_h(\a)$ and $\b\ne\th_h$. Then it is clear 
that $(\th_h-\b,(-\th_h+\d)^\vee)=-1$;  moreover, $\th_h-\b$ is positive 
and belongs, by definition, to $\dperp$. It follows that 
$\th_h-\b\in \ov\Psi_h(\a)$.
\par
By Proposition 7.1 (3) and (5) we have that $|\ov\Psi_h(\a)|=g_h(\a)-1$.
\endemo
\bigskip

\proclaim{Theorem 8.13}
Let $\a\in\Dp_\ell$.
Then $\gom_h(\a)\sbs \gom_{h+1}(\a)$, for
$1\leq h < k(\a)$, and $\gom_{k(\a)}(\a)=\gom(\a)$. Moreover, 
$$\dim\gom_h(\a)=g-1+\sum\limits_{i=1}^{h-1} (g_i(\a)-1).$$
\endproclaim

\demo{Proof}
By Proposition 8.11, $w_{\gom_{h+1}(\a)}=\wa\wah$, for $0\leq
h<k(\a)$; moreover, by Lemma 8.6, $N(\wa\wah)=N(\wa)\cup\wa N(\wah)$.
We have $N(w_{0, \a})=\emptyset$, and $N(\wah)=-\Phi_h(\a)+\d$ for
$1\leq h<k(\a)$; since $\Phi_{h-1}(\a)\sbs \Phi_h(\a)$, we obtain that
$\wa w_{h-1,\a}<\wa\wah$, for $1\leq h<k(\a)$. Therefore
$\gom_h(\a)\sbs \gom_{h+1}(\a)$. It follows that $\gom_{k(\a)}(\a)=
\max \II(\a)$, 
hence, by Corollary 8.8, $\gom_{k(\a)}(\a)=\gom(\a).$
\par
By Lemma 8.6, we have that $\dim\gom_h(\a)=
\ell(\wa)+\ell(w_{h-1,\a})$. By 7.1 (3), $\ell(\wa)=g-1$. Moreover,
$|N(w_{h-1,\a})|=|\Phi_{h-1}(\a)|=\big|\bigcup\limits_{i=1}^{h-1}
\Psi_i\big|$, and since the sets $\Psi_i(\a)$ are pairwise disjoint,
by 8.12, we obtain that $\ell(w_{h-1,\a})=\sum\limits_{i=1}^{h-1}
(g_i(\a)-1)$.
\endemo
\bigskip

\proclaim{Lemma 8.14} Assume $\a\in \Pi_\ell(\t)$, $\t\in X_\a$, 
and $\goj\in\II$.
If there exists $\goi\in \II(\t, \a)$, such that $\goj\spsq \goi$,
then $w_\goi^{-1}w_\goj\in \wap$. In particular $w_\goj(\a)=-\th+2\d$.
\endproclaim

\demo{Proof}
By assumption $w_\goj\geq w_\goi$, in the weak order, hence there
exists $x\in\hatw$ such that $w_\goj=w_\goi x$ and
$\ell(w_\goj)=\ell(w_\goi)+\ell(x)$. Moreover $N(w_\goj)=
N(w_\goi)\cup w_\goi N(x)$, with $w_\goi N(x)\sbsq \Dap$, and indeed,
since $w_\goj\in \wab$, $w_\goi N(x)\sbsq \Dm+\d$. Assume by
contradiction that $x\not\in \wap$ and set $x=x' s_\b
x''$, with $x'\in \wap$, $\b\in\Pia\meno
\a^\perp$, and $\ell(x)=\ell(x')+1+\ell (x'')$, so that $x'(\b)\in
N(x)$. Notice that then $\b\ne\a$, otherwise we have $\a=x'(\a)\in
N(x)$ and hence $-\th+2\d=w_\goi(\a)\in N(w_\goj)$, which is impossible.
Thus, since $\a$ is simple and long, we have $(\b, \a^\vee)=-1$, hence 
 $(w_\goi x'(\b),w_\goi x'(\a))= (w_\goi x'(\b), -\thvee+2\d)=-1$.
But this implies that $w_\goi x'(\b)\in \Dp+\ganz\d$, which is impossible 
since $w_\goi x'(\b)\in N(w_\goj)\sbsq \Dm+\d$. We have therefore that 
$x \in \wap$.
\endemo
\bigskip

Notice that if $\a\in\D_\ell(\t)$ and $\t=\thvee_1+\cdots+\thvee_k$ is 
the standard orthogonal decomposition of $\t$, then $\t=\t_k(\a)$. 
\bigskip
 
\proclaim{Corollary 8.15} 
Assume $\t\in X^+$, $\a,\b\in\Pi_\ell(\t)$, $\a\ne\b$, 
$\goi\in\II(\t,\a)$, and $\goj\in\II(\t,\b)$. Then $\goi,\,\goj$ are not
comparable. In particular, if $\t=\thvee_1+\cdots+\thvee_k$ is 
the standard orthogonal decomposition of $\t$, the maximal elements of 
$\II(\t)$ are exactly the ideals $m_k(\a)$ with $\a\in\Pi_\ell(\t)$.
\endproclaim

\demo{Proof} 
Assume by contradiction $\i\sbsq \goj$. Then by 8.14
$\goj\in\II(\t,\a)$, which is absurd. The last assertion follows from
Proposition 6.7.
\endemo
\smallskip

\remark{\bf Remark} 
It would be interesting to investigate the inclusion relationships
between ideals in $\II(\t,\a),\,\II(\t,\b)$ for $\a,\b\in\Dp_\ell(\t)$.
Related  results can be found  in \cite{P, Corollary 3.3} and in \cite{St2}; we hope
to deal with this problem in a subsequent paper. 
\endremark
\bigskip 

\proclaim{Proposition 8.16} 
If $\a\in \Pi_\ell$, then $\gom(\a)$ is a maximal abelian ideal of $\gob$.
\endproclaim

\demo{Proof}
If $\goi\spsq \gom(\a)$, then by Lemma 8.14 $w_\goi\in \II(\t, \a)$ 
for some $\t\in X_\a$, hence by Theorem 8.13 $\goi\sbsq \gom(\a)$.
\endemo
\bigskip

\ni
{\bf Remark.} Notice that if $\t=\t_{k(\a)}(\a)$, then $\a\not\perp\t$; 
conversely, if $\t\in X_\a$ and $\t=\t_h(\a)$ with $1\leq h<k(\a)$, then
$\a\perp\t$. Hence the maximal ideals of $\gob$ which belong to $\II(\t)$
are exactly the $\gom(\a)$ with $\a\in \Pi_\ell(\t)$ and $\a\not\perp\t$. 
\bigskip

\proclaim{Corollary 8.17} (Panyushev). 
The map $\a\mapsto\gom(\a)$ is a bijection between long simple roots
and maximal abelian ideals of $\frak b$.
\endproclaim

\demo{Proof}
The claim follows from  Proposition 8.16, Corollary 8.15, 
and Theorems 4.4, 8.13.
\endemo
\vskip 1cm

\heading \S9 
Dimension formulas for $\gom_h(\a)$
\endheading
\medskip

In this section we determine the dimension of $\gom_h(\a)$, for all
$\a\in\Pi_\ell$ and $1\leq h\leq k(\a)$.
This calculation has the following
interesting consequence. It is easy  to see that the maximal dimension
$d$ of a maximal commutative subalgebra of $\gog$ coincides with the maximal
dimension of an abelian ideal of $\gob$ (for a proof see \cite{R,
\S7}). Hence
$$d=\mathop{\max}\limits_{\a\in\Pi_\ell} \left\{\dim\gom(\a)\right\}.$$
As expected from Theorem 8.13 we shall see that the maximum is
obtained when $\a$ is a long simple root having maximal  distance from
the affine node in the extended Dynkin diagram of $\D$. The dimensions
$\dim\gom_h(\a)$ will be uniformly calculated using Theorem 8.13:
these calculations afford, as a special case, a conceptual explanation
of the case by case calculations in \cite{PR, Tables I, II}, and in
turn, a proof of Malcev's classical result on the dimension of a
maximal commutative subalgebra \cite{M}. We use for Dynkin diagrams
the numbering of \cite{B}.
\smallskip
  
Before looking at the ideals $\gom_h(\a)$ for each type of root
system, we note that, in general,  for any long simple root
$\a$ we have
$\gom_1(\a)=\max\II(\thvee,\a)$ and \break
$\dim\gom_1(\a)=g-1.$ 
\bigskip

\nd $Type\ A_n$
\smallskip
In this case the element of $X^+$ are linearly ordered and are $\lfloor\frac{n+1}{2}\rfloor$ in number.
Indeed, there exists a unique maximal good chain $\{\D_1=\D\sps\cdots\sps\D_{\lfloor\frac{n+1}{2}\rfloor}\}$; moreover,  any
good chain is an initial section of the maximal one. We have $\D_i=\D(\a_i,\ldots,\a_{n-i+1})\cong A_{2n-i+2}$;
if $\thvee_i$ is the highest root of $\D_i$ and  $\t_i=\thvee_1+\cdots+\thvee_i$
we have  $\t_i=\o_i+\o_{n-i+1}$. Clearly,  $\t_i<\t_j$ when $i<j$. 
If $\a_h\in\D_i$, then  $\widehat\D_k(\a_h)\cong A_{n-2k}$
for $1\leq k \leq i-1$.
Hence for $1\leq i\leq \lfloor\frac{n+1}{2}\rfloor$ and $i\leq h\leq n-i+1$
$$\dim\left(\gom_i(\a_h)\right)=n+\sum_{k=1}^{i-1}(n-2k)=i(n-i+1).$$
We have $d=\lfloor\frac{(n+1)^2}{4}\rfloor$; this value is obtained for
$i=\lfloor\frac{n+1}{2}\rfloor$ and $h=\lfloor\frac{n+1}{2}\rfloor$ if $n$ is odd and  
$h=\frac{n}{2},\,\frac{n}{2}+1$ if $n$ is even.
\bigskip

\nd $Type\ C_n$
\smallskip
Again $X^+$ is linearly ordered and consists of $n$ elements $\t_1=\thvee<\ldots<\t_n=2\o_n$; moreover $\Pi_\ell=\{\a_n\}$
and, for $1\leq i\leq n-1$,  $\widehat\D_i(\a_n)\cong C_{n-i}$, so that 
$$\dim\left(\gom_i(\a_n)\right)=n+\sum_{k=1}^{i-1}(n-k)=\frac{(2n-i+1)i}{2}.$$
We have $d=\frac{n(n+1)}{2}=\dim(\gom_n(\a_n))$.
\bigskip

\nd $Type\ B_n$
\smallskip
In this case there are exactly two maximal good chains. The first is
 $\{\D_1=\D\sps\cdots\sps\D_{\lfloor\frac{n}{2}\rfloor}\}$ with $\D_i=\D(\a_{2i-1},\ldots, \a_n)$ and the other is 
$\{\D_1=\D\sps\tilde \D_2\}$ where $\tilde \D_2$ is the rank one system generated by $\a_1$; this the
good chain of
$\a_1$.
\parni 
The good chain of $\a_2$ consists only of $\D$. 
 We have
$$\align
&\dim\left(\gom_1(\a_h)\right)=2n-2,\quad 1\leq h\leq n-1,\\
&\dim\left(\gom_2(\a_1)\right)=2n-1.\endalign$$
For $h\geq 3$
$$\widehat\D_k(\a_h)\cong\cases A_1\quad&h=2k+1,\\ A_3\quad&h=2k+2,\\ D_{h-2k+1}\quad&2k+3\leq h\leq n-1.\endcases
$$
Hence, for $3\leq h\leq n-1,\,i\geq 2$
$$\dim\left(\gom_i(\a_h))\right)=2n-2+\cases\sum_{k=1}^{i-2}(2h-4k-1)+1\quad&h=2i-1,\\
\sum_{k=1}^{i-2}(2h-4k-1)+3\quad&h=2i,\\
\sum_{k=1}^{i-1}(2h-4k-1)\quad&h>2i-1.\endcases$$
Note that the second summand in the r.h.s. (i.e., $\ell(\check w_{i,\a_h})$) is independent from $n$. The previous calculation
can be expressed as follows
$$\dim(\gom_i(\a_h))=\cases\frac{4n+h^2-3h-2}{2}\quad\!&\text{if the ideal is maximal,}
\\2(n-i^2)+(i-1)(2h+1)\quad&\text{otherwise.}\endcases$$
For $n=3$ we have $d=5$, obtained for $h=1,\,i=2$. For $n\geq 4$ we have
$d=\frac{n(n-1)}{2}+1$ and this value is obtained for $h=n-1$ and $i=\lfloor\frac{n}{2}\rfloor$ (for $n=4$ it is also obtained
for $h=1,\,i=2$).
\smallskip
\nd $Type\ D_n$
\smallskip
The structure of $X^+$ is exactly the same as for $B_n$. The only difference is that any good chain
starts with a system of type $D_n$, hence the only difference in the calculations lies in the Coxeter number of the whole
system. The final results are as follows.
$$\align
&\dim\left(\gom_1(\a_h)\right)=2n-3,\quad 1\leq h\leq n-1,\\
&\dim\left(\gom_2(\a_1)\right)=2n-2,\endalign$$
and for $3\leq h\leq n-1,\,i\geq 2$
$$\dim(\gom_i(\a_h))=\cases\frac{4n+h^2-3h-4}{2} \!\!\!&\text{if the ideal is maximal,}
\\2(n-i^2)+(i-1)(2h+1)-1\!\!\!&\text{otherwise.}\endcases$$

\smallskip
\nd $Type\ E_6$
\smallskip

$X^+$ consists of two elements $\t_1=\thvee<\t_2=\thvee+\thvee(A_5)$ and $\Pi(\t_2)=\Pi\meno \{\a_2\}$.
We have $$\widehat\D_1(\a_i)\cong\cases A_5\quad &i=1,6,\\ A_2 \quad &i=3,5,\\ A_1 \quad &i=4.\endcases$$ therefore
$$\align
&\dim\left(\gom_1(\a_h)\right)=11,\,\quad 1\leq h\leq 6\\
&\dim\left(\gom_2(\a_1)\right)=
\dim\left(\gom_2(\a_6)\right)=16,\\
&\dim\left(\gom_2(\a_3)\right)=
\dim\left(\gom_2(\a_5)\right)=13,\\
&\dim\left(\gom_2(\a_4)\right)=12.\endalign$$
\smallskip
\nd $Type\ E_7$
\smallskip

$X^+$ consists of three elements $\t_1=\thvee<\t_2=\thvee+\thvee(D_6)<\t_3=\thvee+\thvee(D_6)+\thvee(A_1)$ and
$\Pi(\t_2)=\Pi\meno
\{\a_1\},\,\Pi(\t_3)=\{\a_7\}$. We have $$\widehat\D_1(\a_i)\cong\cases A_1\quad &i=2,\\ A_3 \quad &i=3,5,\\ A_2 \quad
&i=4,\\A_5\quad &i=6,\\D_6\quad &i=7,\endcases\qquad\widehat\D_2(\a_7)\cong A_1.$$
therefore
$$\align
&\dim\left(\gom_1(\a_h)\right)=17,\,\quad 1\leq h\leq 7\\
&\dim\left(\gom_2(\a_2)\right)=18,\\
&\dim\left(\gom_2(\a_3)\right)=
\dim\left(\gom_2(\a_5)\right)=20,\\
&\dim\left(\gom_2(\a_4)\right)=19,\\
&\dim\left(\gom_2(\a_6)\right)=22,\\
&\dim\left(\gom_2(\a_7)\right)=26,\\
&\dim\left(\gom_3(\a_7)\right)=27.\endalign$$

\smallskip
\nd $Type\ E_8$
\smallskip

$X^+$ consists of two elements $\t_1=\thvee<\t_2=\thvee+\thvee(E_7)$ and
$\Pi(\t_2)=\Pi\meno
\{\a_8\}$. We have $$\widehat\D_1(\a_i)\cong\cases A_7\quad &i=1,\\ A_5 \quad &i=2,3,\\ A_4 \quad
&i=4,\\A_3\quad &i=5,\\A_2\quad &i=6,\\A_1\quad &i=7,\endcases$$
therefore
$$\align
&\dim\left(\gom_1(\a_h)\right)=29,\,\quad 1\leq h\leq 8\\
&\dim\left(\gom_2(\a_1)\right)=36,\\
&\dim\left(\gom_2(\a_2)\right)=
\dim\left(\gom_2(\a_3)\right)=34,\\
&\dim\left(\gom_2(\a_4)\right)=33,\\
&\dim\left(\gom_2(\a_5)\right)=32,\\
&\dim\left(\gom_2(\a_6)\right)=31,\\
&\dim\left(\gom_2(\a_7)\right)=30.\endalign$$

\smallskip
\nd $Type\ F_4$
\smallskip

$X^+$ consists of two elements $\t_1=\thvee<\t_2=\thvee+\thvee(C_3)$ and
$\Pi(\t_2)=\{\a_2,\a_3,\a_4\}$. Since $\widehat\D(\a_2)\cong A_1$ we have
$$\align
&\dim\left(\gom_1(\a_h)\right)=8,\,\quad h=1,2\\
&\dim\left(\gom_2(\a_2)\right)=9.\endalign$$

\smallskip
\nd $Type\ G_2$
\smallskip
We have $X^+=\{\thvee\}$. There  is only one maximal ideal corresponding to the unique simple long root; it has  dimension $3$.

\bigskip
\remark{\bf Remark} It is also possible to obtain the previous results
using  Proposition 7.7. Indeed
$$\dim(\max\II(\a,\t))=\ell(v_\t)+\ell(y_{\t, \a})+\ell(w_0^\a)-\ell(w_0^{\a,\t})$$
where $w_0^\a$,  is  the longest 
elements in $\wp$ and $w_0^{\a,\t}$ is the longest element in $W((\pi\cap\a^\perp)\meno P_{\t,\a})$.\par
Now $\ell(v_\t)=|\D_\t^2|=|\D_\t^{even}|-|\D_\t^0|$, $\ell(\y)= g(\D(\t))-2$ and the last two summands in the r.h.s. of the
previous formula can be easily calculated determining the number of positive roots of the corresponding subsystems.
\par
As an example, we will write down the computations for types $A_n$ and $E_6$ for the maximal ideals,  displaying the results as
follows. Set $\t_\a=\t_{k(\a)}(\a)$.

$$\alignat5
&\text{long simple root
$\a$}\qquad\qquad&&\tau_{\a}\qquad&&\D_{\t_{\a}}^{even}\qquad&&\D(\t_{\a})\qquad\qquad&
&\ell(w_0^\a)\qquad\\ &dim(\i_\a) && &&\D_{\t_\a}^{0}
&& &&\ell(w_0^{\a,\t_\a})\endalignat$$
 We use  the following  notational conventions: if in the rightmost column appears a root system
of negative rank, then the whole column must not be considered (e.g., type $A_n,\,k=1$). Moreover  a root system of rank zero
should be regarded as the empty set.
\bigskip

\nd $Type\ A_n$
\smallskip

For $1\leq k\leq n$ we have
$$\alignat5
&\a_k\quad&&\o_k+\o_{n-k+1}\quad&& A_{2k-1}\times
A_{n-2k}\qquad&&A_{n-2k+2}\qquad && A_{n-k-1}\\ &  k(n-k+1)&&
&&A_{k-1}\times A_{n-2k}\times A_{k-1} && &&A_{n-2k}\times
A_{k-2}\endalignat$$ \smallskip

\smallskip
\nd $Type\ E_6$
\smallskip
 $$\alignat5
&\a_1,\,\a_6\qquad&&\o_1+\o_6\qquad\quad&&D_5\qquad&&A_5\qquad &&
A_4\\ &16 && &&D_4 && &&A_3\\
&\a_3,\,\a_5\qquad&&\o_1+\o_6\qquad\quad&&D_5\qquad&&A_5\qquad &&
A_1\times A_2\\ &13 && &&D_4 && &&A_2\\
&\a_4\qquad&&\o_1+\o_6\qquad\quad&&D_5\qquad&&A_5\\ &12 && &&D_4
&&\\ &\a_2\qquad&&\o_2\qquad\quad&& A_1\times A_5\qquad&&E_6\\ &11
&& &&A_5
\endalignat$$ 
\endremark

\medskip
\heading{Acknowledgment}\endheading
\bigskip\ni
The authors wish to thank  Corrado De Concini for several useful discussions.

\enddocument
\bye